\documentclass[acmtoms]{acmtrans2m}

\usepackage{graphicx}
\usepackage{psfrag}
\usepackage{url}
\usepackage{amsmath}
\usepackage{amssymb}
\usepackage{algorithm}

\newcommand{\comp}[1]{{\scriptstyle [#1]}}
\newcommand{\dx}{\, \mathrm{d}x}
\newcommand{\dX}{\, \mathrm{d}X}
\newcommand{\R}{\mathbb{R}}
\newcommand{\tab}{\hspace*{2em}}

\title{A Compiler for Variational Forms}

\author{Robert C. Kirby \\ The University of Chicago
        \and
	Anders Logg \\ Toyota Technological Institute at Chicago}

\begin{abstract}
  As a key step towards a complete automation of the finite element
  method, we present a new algorithm for automatic and efficient
  evaluation of multilinear variational forms. The algorithm has been
  implemented in the form of a compiler, the FEniCS Form Compiler
  FFC. We present benchmark results for a series of standard
  variational forms, including the incompressible Navier--Stokes
  equations and linear elasticity. The speedup compared to the
  standard quadrature-based approach is impressive; in some cases the
  speedup is as large as a factor 1000.
\end{abstract}




\begin{document}

\begin{bottomstuff}
Robert C. Kirby, Department of Computer Science, University of Chicago,
1100 East 58th Street, Chicago, Illinois 60637, USA.
\emph{Email:} \texttt{kirby@cs.uchicago.edu}.
This work was supported by the United States Department of Energy
under grant DE-FG02-04ER25650.
\newline
Anders Logg, Toyota Technological Institute at Chicago, University Press Building,
1427 East 60th Street, Chicago, Illinois 60637, USA.
\emph{Email:} \texttt{logg@tti-c.org}.
\end{bottomstuff}

\maketitle

\section{Introduction}

The finite element method provides a general mathematical framework
for the solution of differential equations and can be viewed as a
machine that automates the discretization of differential equations;
given the variational formulation of a differential equation, the
finite element method generates a discrete system of equations for the
approximate solution.

This generality of the finite element method is seldom reflected
in codes, which are often very specialized and can only solve
one particular differential equation or a small set of differential
equations.

There are two major reasons that the finite element method has yet to
be fully automated; the first is the complexity of the task itself,
and the second is that specialized codes often outperform general
codes.  We address both these concerns in this paper.

A basic task of the finite element method is the computation of the
element stiffness matrix from a bilinear form on a local element. In
many applications, computation of element stiffness matrices accounts
for a substantial part of the total run-time of the code.
\cite[SISC]{logg:article:07} This routine is a small amount of code, but it
can be tedious to get it both correct and efficient. While the
standard quadrature-based approach to computing the element stiffness
matrix works on very general variational forms, it is well-known that
precomputing certain quantities in multilinear forms can improve the
efficiency of building finite element matrices.

The methods discussed in this paper for efficient computation of
element stiffness matrices are based on ideas previously presented in
\cite[SISC]{logg:article:07} and \cite[BIT]{KirKne05}, where the basic
idea is to represent the element stiffness matrix as a tensor product.
A similar approach has been implemented earlier in the finite element
library DOLFIN~\cite{logg:www:01,logg:preprint:06}, but only for
linear elements.  The current paper generalizes and formalizes these
ideas and presents an algorithm for generation of the tensor
representation of element stiffness matrices and for evaluation of the
tensor product.  This algorithm has been implemented in the form of
the compiler FFC for variational forms; the compiler takes as input a
variational form in mathematical notation and automatically generates
efficient code (C or C++) for computation of element stiffness
matrices and their insertion into a global sparse matrix. This
includes the generation of code both for the computation of element
stiffness matrices and local-to-global mappings of degrees of freedom.

\subsection{FEniCS and the Automation of CMM}

\begin{figure}[htbp]
  \begin{center}
    \psfrag{in1}{$R(u) = 0$}
    \psfrag{in2}{$\mathrm{TOL} > 0$}
    \psfrag{ut}{$U \approx u$}
    \includegraphics[width=10cm]{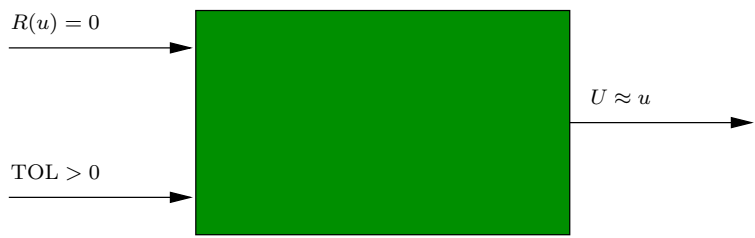}
    \caption{The Automation of Computational Mathematical Modeling.}
    \label{fig:automation}
  \end{center}
\end{figure}

FFC, the FEniCS Form Compiler, is a central component of FEniCS
\cite{logg:www:03}, a project for the Automation of Computational
Mathematical Modeling (ACMM). The central task of ACMM, as formulated
in \cite{logg:thesis:03}, is to create a machine that takes as input a
model $R(u) = A(u) - f$, a tolerance $\mathrm{TOL} > 0$ and a norm $\|\cdot\|$
(or some other measure of quality), and produces as output an
approximate solution $U \approx u$ that satisfies the accuracy
requirement $\| U - u \| < \mathrm{TOL}$ using a minimal amount of
work (see Figure \ref{fig:automation}).  This includes an aspect of
\emph{reliability} (the produced solution should satisfy the accuracy
requirement) and an aspect of \emph{efficiency} (the solution should
be obtained with minimal work).

In many applications, several competing models are under
consideration, and one would like to computationally compare them.
Developing separate, special purpose codes for each model is
prohibitive.  Hence, a key step of ACMM is the \emph{automation of
discretization}, i.e., the automatic translation of a differential
equation into a discrete system of equations, and as noted above this
key step is automated by the finite element method. The FEniCS Form
Compiler FFC may then be viewed as an important step towards the
automation of the finite element method, and thus as an important step
towards the Automation of CMM.

FEniCS software is free software. In particular, FFC is licensed under
the GNU General Public License \cite{www:GPL}. All FEniCS software is
available for download on the FEniCS web site \cite{logg:www:03}.  In
Section \ref{sec:toolchain}, we return to a discussion of the
different components of FEniCS and their relation to FFC.

\subsection{Current finite element software}

Several emerging projects seek to automate important aspects of the
finite element method. By developing libraries in existing languages
or new domain-specific languages, software tools may be built that
allow programmers to define variational forms and other parts of a
finite element method with succinct, mathematical syntax. Existing C++
libraries for finite elements include
DOLFIN~\cite{logg:www:01,logg:preprint:06}, Sundance~\cite{Lon03},
deal.II~\cite{www:deal}, Diffpack~\cite{Lan99} and FEMSTER
\cite{CasKon04,CasRie05}. Projects developing domain-specific
languages for finite element computation include
FreeFEM~\cite{www:FreeFEM} and GetDP~\cite{www:GetDP}. A precursor to
the FEniCS project, Analysa~\cite{www:Analysa}, was a Scheme-like
language for finite element methods. Earlier work on object-oriented
frameworks for finite element computation include \cite{Mac92} and
\cite{MasUsm97}.

While these tools are effective at exploiting modern software
engineering to produce workable systems, we believe that additional
mathematical insight will lead to even more powerful codes with more
general approximating spaces and more powerful algorithms. The FEniCS
project is more ambitious than to just collect and implement existing
ideas.

\subsection{Design goals}

The primary design goal for FFC is to accept as input ``any''
multilinear variational form and ``any'' finite element, and to
generate code that will run with close to optimal performance.

We will make precise below in Section \ref{sec:representation} which
forms and which elements the compiler can currently handle (general
multilinear variational forms with coefficients over affine simplices).

A secondary goal for FFC is to create a new standard in form
evaluation; hopefully FFC can become a standard tool for practitioners
solving partial differential equations using the finite element
method. In addition to generating very efficient code for evaluation
of the element stiffness matrix, FFC thus takes away the burden of
having to implement the code from the developer. Furthermore, if the
code for the element stiffness matrix is generated by a compiler that
is trusted and has gone through rigorous testing, it is easier to
achieve correctness of a simulation code.

\begin{figure}[htbp]
  \begin{center}
    \psfrag{FFC}{\ \large\textbf{\textsf{FFC}}}
    \psfrag{form}{\hspace{0.55cm}$a(v, u) = \int_{\Omega} \nabla v(x)
    \cdot \nabla u(x) \dx$} \includegraphics[width=10cm]{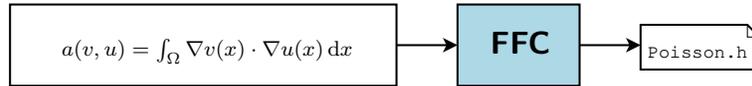}
    \caption{The form compiler FFC takes as input a variational form and
      generates code for evaluation of the form.}
    \label{fig:ffc}
  \end{center}
\end{figure}

The primary output target of FFC is the C++ library DOLFIN. By
default, FFC accepts as input a variational form and generates code
for the evaluation of the variational form in DOLFIN, as illustrated
in Figure \ref{fig:ffc}. Although FFC works closely with other FEniCS
components, such as DOLFIN, it has an abstraction layer that allows
it to be hooked up to multiple backends. One example of this is the
newly added ASE (ANL SIDL Environment, \cite{www:ASE}) format added to FFC,
allowing forms to be compiled for the next generation of PETSc \cite{www:PETSc}.

\subsection{The compiler approach}

It is widely accepted in developing software for
scientific computing that there is a trade-off between generality and
efficiency; a software component that is general in nature, i.e., it
accepts a wide range of inputs, is often less efficient than another
software component that performs the same job on a more limited set of
inputs. As a result, most codes used by practitioners for the solution
of differential equations are very specific.

However, by using a compiler approach, it is possible to combine
generality and efficiency without loss of generality and without loss
of efficiency. This is possible since our compiler works on a very
small family of inputs (multilinear variational forms) with sharply
defined mathematical properties.  Our domain-specific knowledge allows
us to generate much better code than if we used general-purpose
compiler techniques.

\subsection{Outline of this paper}

Before presenting the main algorithm, we give a short background
on the implementation of the finite element method and the evaluation
of variational forms in Section \ref{sec:background}. The main
algorithm is then outlined in Section \ref{sec:evaluation}. In Section
\ref{sec:complexity}, we compare the complexity of form evaluation
for the algorithm used by FFC with the standard quadrature-based
approach. We then discuss the implementation of the form compiler in
some detail in Section \ref{sec:implementation}.

\begin{figure}[htbp]
  \begin{center}
    \includegraphics[width=11cm]{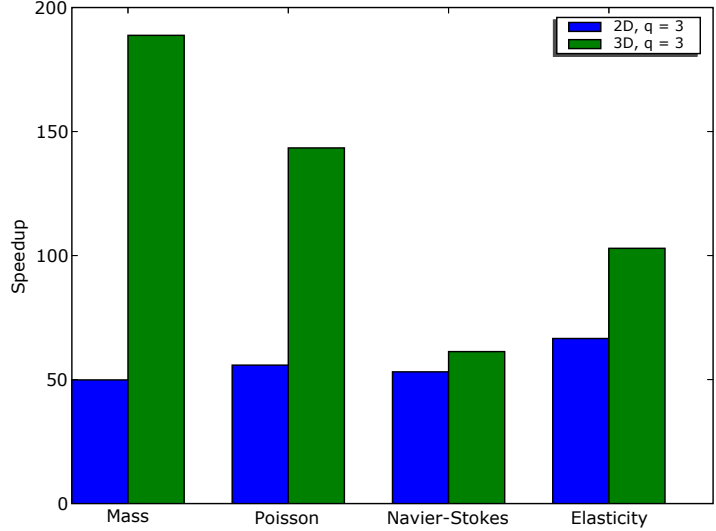}
    \caption{Speedup for a series of standard variational forms (here
      compiled for cubic Lagrange elements on triangles and tetrahedra, respectively).}
    \label{fig:benchmark_3d}
  \end{center}
\end{figure}

In Section \ref{sec:benchmarks}, we compare the CPU time for
evaluating a series of standard variational forms using code
automatically generated by FFC and hand-coded quadrature-based
implementations. The speedup is in all cases significant, in the case
of cubic Lagrange elements on tetrahedra a factor 100 (Figure
\ref{fig:benchmark_3d}).

Finally, in
Section \ref{sec:conclusion}, we summarize the current features and
shortcomings of FFC and give directions for future development and
research.

\section{Background}
\label{sec:background}

In this section, we present a quick background on the finite element
method. The material is standard \cite{Cia76,Hug87,BreSco94,EriEst96},
but is included here to give a context for the presentation of the
form compiler and to summarize the notation used throughout the
remainder of this paper.

For simplicity, we consider here only linear partial differential
equations and note that these play an important role in the
discretization of nonlinear partial differential equations (in Newton
or fixed-point iterations).

\subsection{Variational forms}
We work with the standard variational formulation
of a partial differential equation:
Find $u \in V$ such that
\begin{equation} \label{eq:varform}
  a(v, u) = L(v) \quad \forall v \in \hat{V},
\end{equation}
with $a : \hat{V} \times V \rightarrow \R$ a bilinear form,
$L : \hat{V} \rightarrow \R$ a linear form, and $(\hat{V},V)$ a
pair of suitable function spaces. For the standard example, Poisson's equation
$-\Delta u(x) = f(x)$
with homogeneous Dirichlet conditions on a domain $\Omega$, the bilinear form $a$ is given
by $a(v, u) = \int_{\Omega} \nabla v(x) \cdot \nabla u(x) \dx$, the
linear form $L$ is given by $L(v) = \int_{\Omega} v(x) f(x) \dx$,
and $\hat{V} = V = H^1_0(\Omega)$.

The finite element method discretizes (\ref{eq:varform}) by replacing
$(\hat{V}, V)$ with a pair of (piecewise polynomial) discrete function
spaces. With $\{\hat{\phi}_i\}_{i=1}^M$ a basis for the test space
$\hat{V}$ and $\{\phi_i\}_{i=1}^M$ a basis for the trial space $V$,
we can expand the approximate solution $U$ of (\ref{eq:varform}) in the
basis functions of $V$, $U = \sum_{j=1}^M \xi_j \phi_j$, and obtain
a linear system $A \xi = b$ for the \emph{degrees of freedom}
$\{\xi_j\}_{j=1}^M$ of the approximate solution $U$. The entries of
the matrix $A$ and the vector $b$ defining the linear system are given
by

\begin{equation}
  \begin{split}
    A_{ij} &= a(\hat{\phi}_i, \phi_j), \quad i, j = 1,\ldots,M, \\
    b_i &= L(\hat{\phi}_i), \quad i = 1,\ldots,M.
  \end{split}
\end{equation}

\subsection{Assembly}

The standard algorithm for computing the matrix $A$ (or the vector
$b$) is \emph{assembly}; the matrix is computed by iteration over the
elements $K$ of a triangulation $\mathcal{T}$ of $\Omega$, and the
contribution from each local element is added to the global matrix
$A$.

To see this, we note that if the bilinear form $a$ is expressed as an
integral over the domain $\Omega$, we can write the bilinear form as a
sum of element bilinear forms, $a(v, u) = \sum_{K\in\mathcal{T}} a_K(v,
u)$, and thus
\begin{equation}
  A_{ij} = \sum_{K\in\mathcal{T}} a_K(\hat{\phi}_i, \phi_j),
  \quad i, j = 1,\ldots,M.
\end{equation}
In the case of Poisson's equation, the element bilinear form $a_K$ is
defined by $a_K(v, u) = \int_{K} \nabla v(x) \cdot \nabla u(x) \dx$.

Let now $\{\hat{\phi}_i^K\}_{i=1}^n$ be
the restriction to $K$ of the subset of $\{\hat{\phi}_i\}_{i=1}^M$
supported on $K$ and $\{\phi_i^K\}_{i=1}^n$ the corresponding local basis for
$V$. Furthermore, let $\hat{\iota}(\cdot, \cdot)$ be a mapping from the local numbering scheme to the
global numbering scheme (local-to-global mapping) for the basis functions of $\hat{V}$, so that
$\hat{\phi}_i^K$ is the restriction to $K$ of $\phi_{\hat{\iota}(K,i)}$,
and let $\iota(K, \cdot)$ be the corresponding mapping for $V$.
We may now express an algorithm for the computation of the
matrix~$A$ (Algorithm \ref{alg:assembly,1}).

\begin{algorithm}
  \begin{tabbing}
    $A = 0$\\
    \textbf{for}  {$K \in \mathcal{T}$}\\
    \tab \textbf{for} $i = 1,\ldots,n$ \\
    \tab \tab \textbf{for} {$j = 1,\ldots,n$} \\
    \tab \tab \tab {$A_{\hat{\iota}(K, i)\iota(K,j)} = A_{\hat{\iota}(K,i)\iota(K,j)} + a_K(\hat{\phi}^K_i, \phi_j^K)$} \\
    \tab \tab \textbf{end for} \\
    \tab \textbf{end for} \\
    \textbf{end for}
  \end{tabbing}
  \caption{$A$ = Assemble($a$, $\mathcal{T}$, $\hat{V}$, $V$)}
  \label{alg:assembly,1}
\end{algorithm}

Alternatively, one may define the \emph{element matrix} $A^K$ by
\begin{equation} \label{eq:elementmatrix}
  A^K_{ij} =  a_K(\hat{\phi}^K_i, \phi_j^K)
 \quad i, j = 1,\ldots,n,
\end{equation}
and separate the computation on each element $K$ into two steps:
computation of the element matrix $A^K$ and insertion of $A^K$ into $A$
(Algorithm \ref{alg:assembly,2}).

\begin{algorithm}
  \begin{tabbing}
    {$A = 0$} \\
    \textbf{for} {$K \in \mathcal{T}$} \\
    \tab Compute $A^K$ according to (\ref{eq:elementmatrix}) \\
    \tab Add $A^K$ to $A$ using the local-to-global mappings $(\hat{\iota}(K,\cdot), \iota(K,\cdot))$ \\
    \textbf{end for}
  \end{tabbing}
  \caption{$A$ = Assemble($a$, $\mathcal{T}$, $\hat{V}$, $V$)}
  \label{alg:assembly,2}
\end{algorithm}

Separating the two concerns of computing the element matrix $A^K$ and
adding it to the global matrix $A$ as in Algorithm \ref{alg:assembly,2}
has the advantage that one may
use an optimized library routine for adding the element
matrix $A^K$ to the global matrix $A$.  Sparse matrix libraries such
as PETSc \cite{www:PETSc,BalBus04,BalEij97} often provide optimized routines for this type
of operation. Note that the cost of adding $A^K$ to $A$ may be
substantial even with an efficient implementation of the sparse data
structure for $A$ \cite[SISC]{logg:article:07}.

As we shall see below, we may also take advantage of the separation of
concerns of Algorithm \ref{alg:assembly,2} to optimize the
computation of the element matrix $A^K$. This step is automated by the
form compiler FFC. Given a bilinear (or multilinear) form~$a$, FFC
automatically generates code for run-time computation of the element
matrix $A^K$.

\section{Evaluation of multilinear forms}
\label{sec:evaluation}

In this section, we present the algorithm used by FFC to automatically
generate efficient code for run-time computation of the element matrix $A^K$.

\subsection{Multilinear forms}

Let $\{V_i\}_{i=1}^r$ be a given set of discrete function spaces defined
on a triangulation $\mathcal{T}$ of $\Omega \subset \R^d$.
We consider a general multilinear form $a$ defined on the product space
$V_1 \times V_2 \times \cdots \times V_r$:
\begin{equation}
  a : V_1 \times V_2 \times \cdots \times V_r \rightarrow \R.
\end{equation}
Typically, $r = 1$ (linear form) or $r = 2$ (bilinear form), but
the form compiler FFC can handle multilinear forms of arbitrary \emph{arity} $r$.
Forms of higher arity appear frequently in applications and include
variable coefficient diffusion and advection of momentum in the
incompressible Navier--Stokes equations.

Let now
$\{\phi_i^1\}_{i=1}^{M_1},
 \{\phi_i^2\}_{i=1}^{M_2}, \ldots,
 \{\phi_i^r\}_{i=1}^{M_r}$
be bases of $V_1, V_2, \ldots, V_r$ and let $i = (i_1, i_2, \ldots,
i_r)$ be a multiindex. The multilinear form $a$ then
defines a rank $r$ tensor given by
\begin{equation}
  A_i = a(\phi_{i_1}^1, \phi_{i_2}^2, \ldots, \phi_{i_r}^r).
\end{equation}
In the case of a bilinear form, the tensor $A$ is a matrix (the
stiffness matrix), and in the case of a linear form, the tensor $A$ is
a vector (the load vector).

As discussed in the previous section, to compute the tensor $A$ by
assembly, we need to compute the \emph{element tensor} $A^K$ on each
element $K$ of the triangulation $\mathcal{T}$ of $\Omega$. Let
$\{\phi^{K,1}_i\}_{i=1}^{n_1}$ be the restriction to $K$ of the subset
of $\{\phi_i^1\}_{i=1}^{M_1}$ supported on $K$ and define the local
bases on $K$ for $V_2, \ldots, V_r$ similarly. The rank $r$ element
tensor $A^K$ is then defined by
\begin{equation}
  A^K_i = a_K(\phi_{i_1}^{K,1}, \phi_{i_2}^{K,2}, \ldots, \phi_{i_r}^{K,r}).
\end{equation}

\subsection{Evaluation by tensor representation}
\label{sec:representation}

The element tensor $A^K$ can be efficiently computed by representing
$A^K$ as a special tensor product. Under some mild assumptions which
we shall make precise below, the element tensor $A^K$ can be
represented as the tensor product
of a \emph{reference tensor} $A^0$ and a \emph{geometry tensor} $G_K$:
\begin{equation} \label{eq:tensorproduct}
  A^K_{i} = A^0_{i\alpha} G_K^{\alpha},
\end{equation}
or more generally a sum $A^K_{i} = A^{0,k}_{i\alpha} G_{K,k}^{\alpha}$
of such tensor products, where $i$ and $\alpha$ are multiindices and
we use the convention that repetition of an index means summation over
that index. The rank of the reference tensor is the sum of the rank $r
= |i|$ of the element tensor and the rank $|\alpha|$ of the geometry
tensor $G_K$. As we shall see, the rank of the geometry tensor
depends on the specific form and is typically a function of the
Jacobian of the mapping from the reference element and any
variable coefficients appearing in the form.

Our goal is to develop an algorithm that converts an abstract
representation of a multilinear form into (i)~the values of the
reference tensor~$A^0$ and (ii)~an expression for evaluating the
geometry tensor~$G_K$ for any given element $K$. Note that $A^0$ is
fixed and independent of the element $K$ and may thus be
precomputed. Only $G_K$ has to be computed on each element. As we
shall see below in Section \ref{sec:complexity}, for a wide range of
multilinear forms, this allows for computation of the element tensor
$A^K$ using far fewer floating-point operations than if the element
tensor $A^K$ were computed by quadrature on each element $K$.

To see how to obtain the tensor representation
(\ref{eq:tensorproduct}), we fix a small set of operations, allowing
only multilinear forms that can be expressed through these operations,
and observe how the tensor representation (\ref{eq:tensorproduct})
transforms under these operations. As we shall see below, this covers
a wide range of multilinear forms (but not all). We first develop the
tensor representation in abstract form and then present a number of
test cases that exemplify the general notation in Section
\ref{sec:examples}.

As basic elements, we take the local basis functions
$\{\phi_{\gamma}\}_{\gamma} = \cup_i \{\phi^{K,i}_j\}_{j=1}^{n_i}$
for a set of finite element spaces $V_i$, $i = 1,2,\ldots$, including
the finite element spaces $V_1, V_2, \ldots, V_r$ on which the
multilinear form is defined. Allowing addition $\phi_1 + \phi_2$
and multiplication with scalars $\alpha \phi$,
we obtain a vector space $\mathcal{A}$ of linear combinations of basis functions.
Since $\phi_1 - \phi_2 =
\phi_1 + (-1) \phi_2$ and $\phi/\alpha = (1/\alpha) \phi$,
we can easily equip the vector space with subtraction and division by
scalars.

We next equip our vector space with multiplication between elements
of the vector space. We thus obtain an algebra
$\mathcal{A}$ of linear combinations of products of basis functions.
Finally, we extend $\mathcal{A}$ by adding differentiation $\partial /
\partial x_i$ with respect to a coordinate direction $x_i$,
$i=1,\ldots,d$, on~$K$, to obtain
\begin{equation}
  \mathcal{A} = \{ v : v = \sum
  c_{(\cdot)} \prod
  \frac{\partial^{|(\cdot)|} \phi_{(\cdot)}}
       {\partial x_{(\cdot)}} \},
\end{equation}
where $(\cdot)$ represents some multiindex.

To summarize, $\mathcal{A}$ is the algebra of linear combinations of
products of basis functions or derivatives of basis functions that is
generated from the set of basis functions through addition ($+$),
subtraction ($-$), multiplication $(\cdot)$, including multiplication
with scalars, division by scalars $(/)$, and differentiation
$\partial / \partial x_i$. Note that if the basis
functions are vector-valued (or tensor-valued), the algebra is generated
from the set of scalar components of the basis functions.

We may now state precisely the multilinear forms that the form
compiler FFC can handle, namely those multilinear forms that can be
expressed as integrals over $K$ (or the boundary of $K$) of elements
of the algebra $\mathcal{A}$. Note that not all integrals over $K$ of
elements of $\mathcal{A}$ are multilinear forms; in particular, each
product needs to be linear in each argument of the form.

\begin{figure}[htbp]
  \begin{center}
    \psfrag{p0}{$X^1 = (0,0)$}
    \psfrag{p1}{$X^2 = (1,0)$}
    \psfrag{p2}{$X^3 = (0,1)$}
    \psfrag{xi}{$X$}
    \psfrag{x}{$x = F_K(X)$}
    \psfrag{F=}{}
    \psfrag{F}{$F_K$}
    \psfrag{x0}{$x^1$}
    \psfrag{x1}{$x^2$}
    \psfrag{x2}{$x^3$}
    \psfrag{K0}{$K_0$}
    \psfrag{K}{$K$}
    \includegraphics[width=11cm]{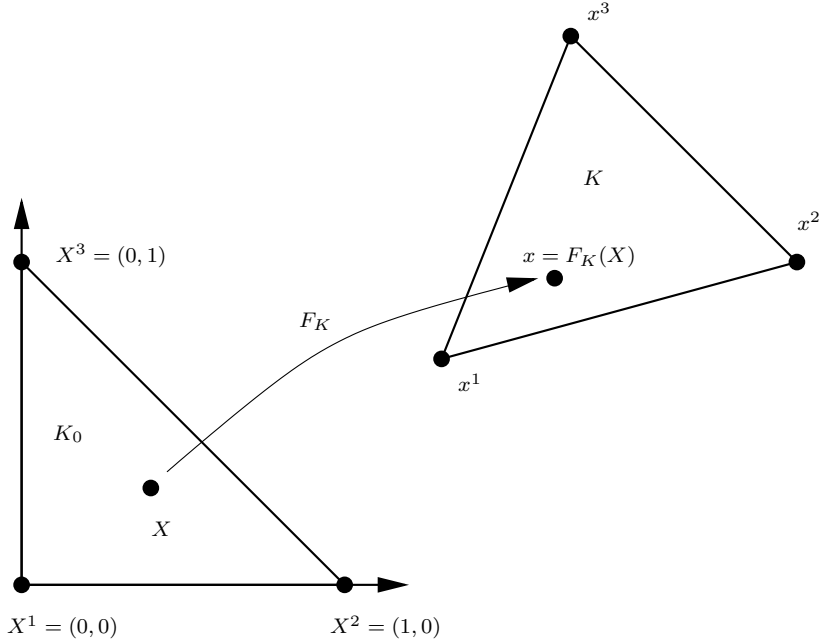}
    \caption{The affine mapping $F_K$ from the reference element $K_0$ to
      the current element $K$.}
    \label{fig:affinemap}
  \end{center}
\end{figure}

The tensor representation (\ref{eq:tensorproduct}) now follows by
a standard change of variables using an affine mapping $F_K : K_0
\rightarrow K$ from a reference element $K_0$ to the current element
$K$ (see Figure \ref{fig:affinemap}). With $\{\Phi_{\gamma}\}_{\gamma}$ the
basis functions on the reference element corresponding to $\{\phi_{\gamma}\}_{\gamma}$, defined
by $\Phi_{\gamma} = \phi_{\gamma} \circ F_K$, we obtain the following
representation of the element tensor $A^K$ corresponding to
$v_i = \left( \sum c_{(\cdot)} \prod
\frac{\partial^{|(\cdot)|} \phi_{(\cdot)}}{\partial x_{(\cdot)}} \right)_i$:
\begin{equation}
  \begin{split}
    A^K_i
    &= a_K(\phi_{i_1}^{K,1}, \phi_{i_2}^{K,2}, \ldots, \phi_{i_r}^{K,r})
    = \int_K v_i \dx \\
    &= \left( \int_K \sum c_{(\cdot)} \prod
    \frac{\partial^{|(\cdot)|} \phi_{(\cdot)}}{\partial x_{(\cdot)}} \dx \right)_i
    = \sum \left( c_{(\cdot)} \int_K \prod
    \frac{\partial^{|(\cdot)|} \phi_{(\cdot)}}{\partial x_{(\cdot)}} \dx \right)_i \\
    &= \sum \left( c_{(\cdot)} \det F_K' \prod
    \frac{\partial X_{(\cdot)}}{\partial x_{(\cdot)}} \right)_{\alpha}
    \left( \int_{K_0} \prod
    \frac{\partial^{|(\cdot)|} \Phi_{(\cdot)}}{\partial X_{(\cdot)}} \dX \right)_{i\alpha} \\
    &= A^{0,k}_{i\alpha} G_{K,k}^{\alpha},
  \end{split}
\end{equation}
where
\begin{eqnarray}
  A^{0,k}_{i\alpha}
  &= \left( \int_{K_0} \prod
  \frac{\partial^{|(\cdot)|} \Phi_{(\cdot)}}{\partial X_{(\cdot)}}
  \dX \right)_{i\alpha}, \\
  G_{K,k}^{\alpha}
  &= \left( c_{(\cdot)} \det F_K' \prod
  \frac{\partial X_{(\cdot)}}{\partial x_{(\cdot)}} \right)_{\alpha}.
\end{eqnarray}

We have here used the fact that the mapping is affine to pull the
determinant and transforms of derivatives outside of the integral. For
a discussion of non-affine mappings, including the Piola transform and
isoparametric mapping, see Section~\ref{sec:extensions} below.

Note that the expression for the geometry tensor $G_{K,k}$ implicitly
contains a summation if an index is repeated twice. Also note that
the geometry tensor contains any variable coefficients appearing in
the form.

As we shall see below in Section \ref{sec:implementation}, the
representation of a multilinear form as an integral over $K$ of an element
of $\mathcal{A}$ is automatically available to the form compiler FFC,
since a multilinear form must be specified using the basic operations
that generate $\mathcal{A}$.

\subsection{Test cases}
\label{sec:examples}

To make these ideas concrete, we write down the explicit tensor
representation (\ref{eq:tensorproduct}) of the element tensor $A^K$
for a series of standard forms. We return to these test cases below in
Section~\ref{sec:benchmarks} when we present benchmark results for each test case.

\subsubsection*{Test case 1: the mass matrix}

As a first simple example, we consider the computation of the mass matrix
$M$ with $M_{i_1i_2} = a(\phi^1_{i_1}, \phi^2_{i_2})$ and the bilinear form $a$
given by
\begin{equation} \label{eq:mass}
  a(v, u) = \int_{\Omega} v(x) u(x) \dx.
\end{equation}
By a change of variables given by the affine mapping
$F_K : K_0 \rightarrow K$, we obtain
\begin{equation}
  A^K_i = \int_K \phi^{K,1}_{i_1}(x) \phi^{K,2}_{i_2}(x) \dx
  = \det F_K' \int_{K_0} \Phi^1_{i_1}(X) \Phi^2_{i_2}(X) \dX =
  A^0_{i} G_K,
\end{equation}
where $A^0_{i} = \int_{K_0} \Phi_{i_1}^1(X) \Phi_{i_2}^2(X) \dX$ and
$G_K = \det F_K'$. In this case, the reference tensor $A^0$ is a rank
two tensor (a matrix) and the geometry tensor $G_K$ is a rank zero
tensor (a scalar). By precomputing the reference tensor $A^0$, we may
thus compute the element tensor $A^K$ on each element $K$ by just
multiplying the precomputed reference tensor with the determinant of (the
derivative of) the affine mapping $F_K$.

\subsubsection*{Test case 2: Poisson's equation}

As a second example, consider the bilinear form for Poisson's equation,
\begin{equation} \label{eq:poisson}
  a(v, u) = \int_{\Omega} \nabla v(x) \cdot \nabla u(x) \dx.
\end{equation}
By a change of variables as above, we obtain the following
representation of the element tensor $A^K$:
\begin{equation}
  \begin{split}
    A^K_i &= \int_{K} \nabla \phi^{K,1}_{i_1}(x) \cdot \nabla \phi^{K,2}_{i_2}(x) \dx \\
    &= \det F_K'
    \frac{\partial X_{\alpha_1}}{\partial x_{\beta}}
    \frac{\partial X_{\alpha_2}}{\partial x_{\beta}}
    \int_{K_0}
    \frac{\partial \Phi^1_{i_1}(X)}{\partial X_{\alpha_1}}
    \frac{\partial \Phi^2_{i_2}(X)}{\partial X_{\alpha_2}} \dX
    = A^0_{i\alpha} G_K^{\alpha},
  \end{split}
\end{equation}
where
$A^0_{i\alpha} = \int_{K_0}
\frac{\partial \Phi^1_{i_1}(X)}{\partial X_{\alpha_1}}
\frac{\partial \Phi^2_{i_2}(X)}{\partial X_{\alpha_2}} \dX$
and
$G_K^{\alpha} = \det F_K'
\frac{\partial X_{\alpha_1}}{\partial x_{\beta}}
\frac{\partial X_{\alpha_2}}{\partial x_{\beta}}$.
We see that the reference tensor $A^0$ is here a rank four tensor and
the geometry tensor $G_K$ is a rank two tensor (one index for each
derivative appearing in the form).

\subsubsection*{Test case 3: Navier--Stokes}

Consider now the nonlinear term $u \cdot \nabla u$ of the
incompressible Navier--Stokes equations,
\begin{equation}
  \begin{split}
    \dot{u} + u \cdot \nabla u - \nu \Delta u + \nabla p &= f, \\
    \nabla \cdot u &= 0.
  \end{split}
\end{equation}
To solve the Navier--Stokes equations by fixed-point iteration (see for example
\cite{EriEst03c}), we write the nonlinear term in the form
$u \cdot \nabla u = w \cdot \nabla u$ with $w = u$ and consider $w$ as fixed.
We then obtain the following bilinear form:
\begin{equation}
  a(v, u) = a_w(v, u) = \int_{\Omega} v(x) \cdot (w(x) \cdot \nabla) u(x) \dx.
\end{equation}
Note that we may alternatively think of this as a trilinear form,
$a = a(v, u, w)$.

Let now $\{w^K_{\alpha}\}_{\alpha}$ be the expansion coefficients for $w$ in
a finite element basis on the current element $K$, and let
$v\comp{i}$ denote component $i$ of a vector-valued function~$v$.
Furthermore, assume that $u$ and $w$ are discretized using the
same discrete space $V = V_2$. We then obtain the following
representation of the element tensor $A^K$:
\begin{equation} \label{eq:navierstokes}
  \begin{split}
    A^K_i &= \int_{K} \phi^{K,1}_{i_1}(x) \cdot (w(x) \cdot \nabla) \phi^{K,2}_{i_2}(x) \dx \\
    &= \det F_K' \frac{\partial X_{\alpha_3}}{\partial x_{\alpha_1}} w^K_{\alpha_2}
    \int_{K_0} \Phi^1_{i_1}\comp{\beta}(X) \Phi^2_{\alpha_2}\comp{\alpha_1}(X)
    \frac{\partial \Phi^2_{i_2}\comp{\beta}(X)}{\partial X_{\alpha_3}} \dX
    = A^0_{i\alpha} G_K^{\alpha},
  \end{split}
\end{equation}
where
$A^0_{i\alpha} = \int_{K_0}
\Phi^1_{i_1}\comp{\beta}(X) \Phi^2_{\alpha_2}\comp{\alpha_1}(X)
\frac{\partial \Phi^2_{i_2}\comp{\beta}(X)}{\partial X_{\alpha_3}} \dX$
and
$G_K^{\alpha} = \det F_K'
\frac{\partial X_{\alpha_3}}{\partial x_{\alpha_1}} w^K_{\alpha_2}$.
In this case, the reference tensor $A^0$ is a rank five tensor and the
geometry tensor $G_K$ is a rank three tensor (one index for the
derivative, one for the function $w$, and one for the scalar product).

\subsubsection*{Test case 4: linear elasticity}

Finally, consider the strain-strain term of linear elasticity,
\begin{equation}
  \begin{split}
    a(v, u) &= \int_{\Omega} \frac{1}{4}
    (\nabla v + (\nabla v)^{\top}) : (\nabla u + (\nabla u)^{\top}) \dx \\
    &= \int_{\Omega}
    \frac{1}{4}
    \frac{\partial v_i}{\partial x_j}
    \frac{\partial u_i}{\partial x_j} \dx +
    \frac{1}{4}
    \frac{\partial v_i}{\partial x_j}
    \frac{\partial u_j}{\partial x_i} \dx +
    \frac{1}{4}
    \frac{\partial v_j}{\partial x_i}
    \frac{\partial u_i}{\partial x_j} \dx +
    \frac{1}{4}
    \frac{\partial v_j}{\partial x_i}
    \frac{\partial u_j}{\partial x_i} \dx.
  \end{split}
\end{equation}
Considering here only the first term, a change of
variables leads to the following representation of the element tensor
$A^{K,1}$:
\begin{equation}
  \begin{split}
    A^{K,1}_i &= \int_{K}
    \frac{1}{4}
    \frac{\partial \phi^{K,1}_{i_1}\comp{\beta_1}(x)}{\partial x_{\beta_2}}
    \frac{\partial \phi^{K,2}_{i_2}\comp{\beta_1}(x)}{\partial x_{\beta_2}} \dx \\
    &=
    \frac{1}{4}
    \det F_K'
    \frac{\partial X_{\alpha_1}}{\partial x_{\beta_2}}
    \frac{\partial X_{\alpha_2}}{\partial x_{\beta_2}}
    \int_{K_0}
    \frac{\partial \Phi^1_{i_1}\comp{\beta_1}(X)}{\partial X_{\alpha_1}}
    \frac{\partial \Phi^2_{i_2}\comp{\beta_1}(X)}{\partial X_{\alpha_2}} \dX
    = A^{0,1}_{i\alpha} G_{K,1}^{\alpha},
  \end{split}
\end{equation}
where
$A^{0,1}_{i\alpha} =
\int_{K_0}
\frac{\partial \Phi^1_{i_1}\comp{\beta_1}(X)}{\partial X_{\alpha_1}}
\frac{\partial \Phi^2_{i_2}\comp{\beta_1}(X)}{\partial X_{\alpha_2}} \dX$
and
$G_{K,1}^{\alpha} =
\frac{1}{4}
\det F_K'
\frac{\partial X_{\alpha_1}}{\partial x_{\beta_2}}
\frac{\partial X_{\alpha_2}}{\partial x_{\beta_2}}$.
Here, the reference tensor $A^{0,1}$ is a rank four tensor and the
geometry tensor $G_{K,1}$ is a rank two tensor (one index for each
derivative).

\subsection{Extensions}
\label{sec:extensions}

The current implementation of FFC supports only affinely mapped
Lagrange elements and linear problems, but it is interesting to
consider the generalization of our approach to other kinds of function
spaces such as Raviart-Thomas~\cite{RavTho77a} elements for
$H(\mathrm{div})$ and curvilinear mappings such as arise with
isoparametric elements, as well as how nonlinear problems may also be
automated.

\subsubsection{$H(\mathrm{div})$ and $H(\mathrm{curl})$ conforming elements}

Implementing of $H(\mathrm{div})$ or $H(\mathrm{curl})$ elements requires
two kinds of generalizations to FFC. First of all, the basis functions
are mapped from the reference element by the Piola
transform~\cite{BreFor91} rather than the standard change of
coordinates.  With $F_K : K_0 \rightarrow K$ the standard affine
mapping for an element $K$, $F^\prime_K$ the Fr\'echet derivative of
the mapping and $\det{F^\prime_K}$ its determinant, the Piola mapping
is defined by $ \mathcal{F}_K( \Psi ) =
\frac{1}{\det{F^\prime_K}}F^\prime_K (\Psi \circ \left( F_K
\right)^{-1})$.  Since our tools already track Jacobians and
determinants for differentiating through affine mappings, it should be
straightforward to support the Piola mapping. Second, defining the
mapping between local and global degrees of freedom becomes more
complicated, as we must keep track of directions of vector
components as done in FEMSTER~\cite{CasKon04,CasRie05}.

As an example of using the Piola transform, we consider the
Raviart-Thomas elements with the standard
(vector-valued) nodal basis $\{ \Psi_i\}_{i=1}^d$ on the reference
element. We compute the mass matrix $M$ with $M_{i_1i_2} =
a(\psi_{i_1}, \psi_{i_2})$ and the bilinear form $a$ given by
\begin{equation} \label{eq:masspiola}
  a(v, u) = \int_{\Omega} v(x) \cdot u(x) \dx.
\end{equation}

On $K$, the basis functions are given by $\psi^K_i = \mathcal{F}_K(\Psi_i)$
and computing the element tensor $A^K$, we obtain
\begin{equation}
  \begin{split}
    A^K_i &= \int_{K} \psi^{K}_{i_1}(x) \cdot \psi^{K}_{i_2}(x) \dx \\
    &= \frac{1}{\det F_K'}
    \frac{\partial x_{\beta}}{\partial X_{\alpha_1}}
    \frac{\partial x_{\beta}}{\partial X_{\alpha_2}}
    \int_{K_0}
    \Psi_{i_1}\comp{\alpha_1}\Psi_{i_2}\comp{\alpha_2} \dX
    = A^0_{i\alpha} G_K^{\alpha},
  \end{split}
\end{equation}
where
$A^0_{i\alpha} = \int_{K_0}
    \Psi_{i_1}\comp{\alpha_1}\Psi_{i_2}\comp{\alpha_2}
 \dX$
and
$G_K^{\alpha} = \frac{1}{\det F_K'}
\frac{\partial x_{\beta}}{\partial X_{\alpha_1}}
\frac{\partial x_{\beta}}{\partial X_{\alpha_2}}$.
We see that, like for Poisson, the reference tensor $A^0$ is rank four and
the geometry tensor $G_K$ is rank two.

\subsubsection{Curvilinear elements}

Our techniques may be generalized to cases in which the Jacobian
varies spatially within elements, such as when curvilinear elements or
general quadrilaterals or hexahedra are used. In this case, we can
replace integration with summation over quadrature points and obtain a
formulation based on tensor contraction, albeit with a higher run-time
complexity than for affine elements.

To illustrate this, we consider the case of the very simple
bilinear form
\begin{equation}
  \label{eq:advect}
  a(v, u) = \int_{\Omega} v \frac{\partial u}{\partial x_1} \dx.
\end{equation}

If the mapping from the reference element $K_0$ to an element $K$
is curvilinear, then we will be unable to pull the Jacobian and
derivatives out of the integral. However, we will still be able to
phrase a run-time tensor contraction with an extra index for the
quadrature points. The element matrix for (\ref{eq:advect}) is
\begin{equation}
  A^K_i = \int_{K} \phi_{i_1} \frac{\partial \phi_{i_2}}{\partial x_1} \dx,
\end{equation}
and changing coordinates we obtain
\begin{equation}
A^K_i = \int_{K_0}
\Phi_{i_1}
\frac{\partial \Phi_{i_2}}{\partial X_{\alpha_1}}
\frac{\partial X_{\alpha_1}}{\partial x_1} \det{F_K'}
 \dX.
\end{equation}

Approximating the integral by quadrature, we
let $ \{ X_k \}_{k=1}^N$ be a set of quadrature points
on the reference element $K_0$ with $\{ w_k \}_{k=1}^{N}$ the
corresponding weights. We thus obtain the representation
\begin{equation}
  A^K_{i} \approx \tilde{A}^K_{i} =
  \sum_{k=1}^N
  w_k
  \Phi_{i_1}(X_k)
  \frac{\partial \Phi_{i_2}}{\partial X_{\alpha_1}} (X_k)
  \frac{\partial X_{\alpha_1}}{\partial x_1} (X_k)
  \det{F_K'}(X_k) =
  \tilde{A}^0_{i\alpha} G_K^{\alpha},
\end{equation}
where
\begin{equation}
  \tilde{A}^0_{i\alpha} =
  w_{\alpha_2}
  \Phi_{i_1}(X_{\alpha_2}) \frac{\partial \Phi_{i_2}}{\partial X_{\alpha_1}} (X_{\alpha_2})
\end{equation}
and
\begin{equation}
  G_\alpha =
  \frac{\partial X_{\alpha_1}}{\partial x_1} (X_{\alpha_2}) \det{F_K'}(X_{\alpha_2}).
\end{equation}

Note that $\tilde{A}^0_{\alpha}$ may be computed entirely at compile-time,
as can an expression for $G_\alpha$, whereas the values of
$G_\alpha$ depend on the geometry given at run-time. The tensors to be
contracted at runtime have one extra dimension compared to a situation
where the mapping is affine. Still, this computation (once the
geometry tensor is computed on an element) is readily phrased as a
matrix-vector product. Hence, we have given an example of how the more
traditional way of expressing quadrature and our precomputation are
both instances of a high-rank tensor contraction.  One could interpret
this formulation as saying that quadrature (expressed as we have here)
gives a general model for computation and that precomputation is
possible as a compile-time optimization in the case where the mapping
is affine.

An open interesting problem would be to study under what conditions
one could specialize a system such as FFC to use bases with
tensor-product decompositions (available on unstructured as well as
structured shapes~\cite{KarShe99}) and automatically generate efficient
matrix-vector products as are manually implemented in spectral element
methods.

\subsubsection{Nonlinear forms}

For a nonlinear variational problem,
\begin{equation} \label{eq:varform,nonlinear}
  a(v, u) = L(v) \quad \forall v \in \hat{V},
\end{equation}
with $a$ nonlinear in $u$, we can solve by direct fixed-point
iteration on the unknown $u$ as discussed above for test case~3, or we
can compute the (Fr\'echet) derivative $a'_u$ of the nonlinear form
$a$ and solve by Newton's method. For multilinear forms, defining the
nonlinear residual and constructing the Jacobian can be performed with
the current capabilities of FFC.

To solve the variational problem (\ref{eq:varform,nonlinear}) by
Newton's method, we differentiate with respect to $u$ to obtain a
variational problem for the \emph{increment} $\delta u$:
\begin{equation} \label{eq:varform,newton}
  a'_u(v, \delta u) = - (a(v, u) - L(v)) \quad \forall v \in \hat{V}.
\end{equation}
As an example, consider the nonlinear Poisson's equation $-c(u)\Delta
u = f$ with $c(u) = u$.  For more general $c$, considering the
projection of $c(u)$ into the finite element space leads to a
multilinear form.  Differentiating the form with respect to $u$, we
obtain the following variational problem: Find $\delta u \in V$ such
that
\begin{equation}
  \int_{\Omega} \delta u \nabla v \cdot \nabla u +
  u \nabla v \cdot \nabla \delta u \dx =
  \int_{\Omega} v f \dx - \int_{\Omega} u \nabla v \cdot \nabla u \dx
  \quad \forall v \in \hat{V}.
\end{equation}

Our current capabilities would allow us to define two forms, one to
evaluate the nonlinear residual and another to construct the Jacobian
matrix.  This is sufficient to set up a nonlinear solver. Obviously,
extending FFC to symbolically differentiate the nonlinear form would
be more satisfying. We remark that the code Sundance~\cite{Lon03}
currently has such capabilities. It should be possible to leverage
such tools in the future to combine our precomputation techniques with
automatic differentiation.

\subsection{Optimization}
\label{sec:optimization}

We consider here three different kinds of optimization that could be
built into FFC in the future.  For one, the current code is generated
entirely straightline as a sequence of arithmetic and assignment. It
should be possible to store the tensor $A^0$ in a contiguous array.
Moreover, each $G_K$ may be considered as a tensor or flattened into a
vector. In the latter case, the action of forming the element matrix
for one element may be written as a matrix-vector multiply using the
level 2 BLAS. Once this observation is made, it is straightforward to
see that we could form several $G_K$ vectors and make better use of
cache by computing several element matrices at once by a matrix-matrix
multiply and the level 3 BLAS.

This corresponds to a coarse-grained optimization. In other
work~\cite[SISC]{logg:article:07}, \cite[BIT]{KirKne05}, we have seen
that for many forms, the entries in $A^0$ are related in such ways
that various entries of the element matrices may be formed in fewer
operations.  For example, if two entries of $A^0$ are close together
in Hamming distance, then the contraction of one entry with $G_K$ can
be computed efficiently from the other. As our code for optimization,
FErari, evolves, we will integrate it with FFC as an optimizing
backend. It will be simple to compare the output of FErari to the best
performance using the BLAS, and let FFC output the best of the two
(which may be highly problem-dependent).

Finally, optimizations that arise from the variational form itself
will be fruitful to explore in the future. For example, it should be
possible to detect when a variational form is symmetric within FFC, as
this leads to fewer operations to form the associated matrix.
Moreover, for forms over vector-valued elements that have a Cartesian
product basis (each basis function has support in only one component),
other kinds of optimizations are appropriate. For example, the
viscosity operator for Navier-Stokes is the vector Laplacian, which
can be written as a block-diagonal matrix with one axis for each
spatial dimension. By "taking apart" the basis functions, we hope to
uncover this block structure, which will lead to more efficient
compilation and hopefully more efficient code.

\section{Complexity of form evaluation}
\label{sec:complexity}

We now compare the proposed algorithm based on tensor representation to
the standard quadrature-based approach. As we shall see,
tensor representation can be much more efficient than quadrature for a
wide range of forms.

\subsection{Basic assumptions and notation}

To analyze the complexity of form evaluation, we make the following
simplifying assumptions:
\begin{itemize}
\item
  the form is bilinear, i.e., $r = |i| = 2$;
\item
  the form can be represented as one tensor product, i.e.,
  $A^K_i = A^0_{i\alpha} G_K^{\alpha}$;
\item
  the basis functions are scalar;
\item
  integrals are computed exactly, i.e., the order of the quadrature
  rule must match the polynomial order of the integrand.
\end{itemize}

We shall use the following notation: $q$ is the polynomial order of the
basis functions on every element, $p$ is the total polynomial order of
the integrand of the form, $d$ is the dimension of $\Omega$, $n$ is
the dimension of the function space on an element, and $N$ is the number
of quadrature points needed to integrate polynomials of degree $p$
exactly.

Furthermore, let $n_f$ be the number of functions appearing in the
form. For test cases~1--4 above, $n_f = 0$ in all cases except test
case 3 (Navier--Stokes) where $n_f = 1$. We use $n_D$ to denote the
number of differential operators. For test cases~1--4,
we have $n_D = 0$ in case~1 (the mass matrix), $n_D = 2$ in case~2
(Poisson), $n_D = 1$ in case~3 (Navier--Stokes), and $n_D = 2$ in
case~4 (linear elasticity).

\subsection{Complexity of tensor representation}

The element tensor $A^K$ has $n^2$ entries. The number of basis
functions $n$ for polynomials of degree $q$ in $d$ dimensions is $\sim
q^d$.
To compute each entry $A^K_i$ of the element tensor $A^K$ using
tensor representation, we need to compute the tensor product between
$A^0_{i\cdot}$ and $G_K$. The geometry tensor $G_K$ has rank
$n_f + n_D$ and the number entries of $G_K$ is $n^{n_f} d^{n_D}$.
The cost for computing the $n^2$ entries of the element tensor $A^K$
using tensor representation is thus
\begin{equation}
  T_T \sim n^2 n^{n_f} d^{n_D}
  \sim (q^d)^2 (q^d)^{n_f} d^{n_D}
  \sim q^{2d + n_f d} d^{n_D}.
\end{equation}
Note that there is no run-time cost associated with computing the tensor
representation (\ref{eq:tensorproduct}), since this is computed at
compile-time. Also note that we have not taken into account any of the
optimizations discussed in Section \ref{sec:optimization}. These
optimizations can in some cases significantly reduce the operation
count.

\subsection{Complexity of quadrature}

To compute each entry $A^K_i$ of the element tensor $A^K$ using
quadrature, we
need to evaluate an integrand of total order $p$ at $N$ quadrature
points. The number of quadrature points needed to integrate
polynomials of order $p$ exactly in $d$ dimensions is $N \sim
p^d$. Since the form is bilinear with basis functions of order $q$,
the total order is $p = 2q + n_f q - n_D$. It is difficult to estimate
precisely the cost of evaluating the integrand at each quadrature
point, but a reasonable estimate is $n_f + n_D d + 1$. Note that we
assume that all basis functions and their derivatives have been
pretabulated at all quadrature points on the reference element.

We thus obtain the following estimate of the total cost for computing
the $n^2$ entries of the element tensor $A^K$ using quadrature:
\begin{equation}
  \begin{split}
    T_Q &\sim n^2 N (n_f + n_D d + 1)
    \sim (q^d)^2 p^d (n_f + n_D d + 1) \\
    &\sim q^{2d} (2q + n_f q - n_D)^d (n_f + n_D d + 1).
  \end{split}
\end{equation}

\subsection{Tensor representation vs. Quadrature}

Comparing tensor representation with quadrature, the speedup of tensor
representation is
\begin{equation}
  T_Q / T_T \sim \frac{(2q + n_f q - n_D)^d (n_f + n_D d + 1)}{q^{n_f d} d^{n_D}}.
\end{equation}

We immediately note that there can be a significant speedup for $n_f =
0$, since $T_T/n^2$ is then independent of the polynomial degree
$q$. In particular, we note that for the mass matrix ($n_f = n_D = 0$)
the speedup is $T_Q/T_T \sim (2q)^d$, and for Poisson's equation ($n_f
= 0$, $n_D = 2$) the speedup is $T_Q/T_T \sim (2q - 2)^d (2d + 1) /
d^2$. As we shall see below, the speedup for test cases~1--4 is
significant, even for $q = 0$.

On the other hand, we note that quadrature may be more efficient if
$n_f$ is large. Also, if one takes into account that underintegration
is possible (choosing a smaller $N$ than given by the polynomial
order~$p$), it is less clear which approach is most efficient in any given
case. It is known~\cite{Cia76} that second order elliptic variational
problems using polynomials of degree $q$ require an integration rule
that is exact only on polynomials of degree $2q-2$ to ensure the
proper convergence rate, regardless of the arity of the form. However,
our overall flop count is lower for bilinear and likely trilinear
forms, and at any rate, our code is simpler for compilers to optimize
than quadrature loops.

Ultimately, one may thus imagine an intelligent system that automatically
makes the choice between tensor representation and quadrature in each
specific situation.

\section{Implementation}
\label{sec:implementation}

We now discuss a number of important aspects of the implementation of
the form compiler FFC. We also write down the forms for the test cases
discussed above in Section \ref{sec:examples} in the language of the
form compiler FFC.  Basically, we can consider FFC's functionality
broken into three phases.  First, it takes an expression for a
multilinear form and generates the tensor $A^0$.  While doing
this, it derives an expression for evaluation of the element tensor $G_K$
from the affine mapping and the coefficients of the form. Finally, it generates
code for evaluating $G_K$ and contracting it with $A^0$, and for
constructing the local-to-global mapping.

\subsection{Parsing of forms}
\label{sec:parsing}

The form compiler FFC implements a domain-specific language (DSL) for
variational forms, using Python as the host language. A language of
variational forms is obtained by overloading the appropriate
operators, including addition \texttt{+}, subtraction \texttt{-},
multiplication (\texttt{*}), and differentiation \texttt{.dx($\cdot$)}
for a hierarchy of classes corresponding to the algebra $\mathcal{A}$
discussed above in Section \ref{sec:representation}. FFC thus uses the
built-in parser of Python to process variational forms.

\subsection{Generation of the tensor representation}

The basic elements of the algebra are objects of the class
\texttt{BasisFunction}, representing (derivatives of) basis functions
of some given finite element space. Each \texttt{BasisFunction} is
associated with a particular finite element space and different
\texttt{BasisFunction}s may be associated with different finite
element spaces.  Products of scalars and (derivatives of) basis
functions are represented by the class \texttt{Product}, and sums of
such products are represented by the class \texttt{Sum}. In addition,
we include a class \texttt{Function}, representing linear combinations
of basis functions (coefficients).  In the diagrams of Tables
\ref{tab:unary} and \ref{tab:binary}, we summarize the basic unary and
binary operators respectively implemented for the hierarchy of
classes.

Note that by declaring a common base class for \texttt{BasisFunction},
\texttt{Product}, \texttt{Sum}, and \texttt{Function}, some of the
operations can be grouped together to simplify the implementation. As
a result, most operators will directly yield a \texttt{Sum}. Also note
that the algebra of \texttt{Sum}s is closed under the operations
listed above.

\begin{table}[htbp]
  \begin{center}
    \begin{tabular}{l|llll}
      \texttt{op}           & B & F & P & S \\
      \hline
      \texttt{-}            & P & S & P & S \\
      \texttt{.dx($\cdot$)} & P & S & S & S
    \end{tabular}
    \caption{Unary operators and their results for the
      classes \texttt{BasisFunction} (B), \texttt{Function} (F),
      \texttt{Product} (P), and \texttt{Sum} (S).}
    \label{tab:unary}
  \end{center}
\end{table}

\begin{table}[htbp]
  \begin{center}
    \begin{tabular}{l|llll}
      \texttt{+/-}          & B & F & P & S \\
      \hline
      B                     & S & S & S & S \\
      F                     & S & S & S & S \\
      P                     & S & S & S & S \\
      S                     & S & S & S & S \\
    \end{tabular}
    \hspace{1cm}
    \begin{tabular}{l|llll}
      \texttt{*}            & B & F & P & S \\
      \hline
      B                     & P & S & P & S  \\
      F                     & S & S & S & S  \\
      P                     & P & S & P & S  \\
      S                     & S & S & S & S  \\
    \end{tabular}
    \caption{Binary operators and their results for the
      classes \texttt{BasisFunction} (B), \texttt{Function} (F),
      \texttt{Product} (P), and \texttt{Sum} (S).}
    \label{tab:binary}
  \end{center}
\end{table}

By associating with each object one or more \emph{indices},
implemented by the class \texttt{Index}, an object of type
\texttt{Sum} automatically represents a tensor, and by differentiating
between different types of indices, an object of type \texttt{Sum}
automatically encodes the tensor representation
(\ref{eq:tensorproduct}). FFC differentiates between four different
types of indices: \emph{primary}, \emph{secondary}, \emph{auxiliary},
and \emph{fixed}. A primary index ($i$) is associated with the multiindex of
the element tensor $A^K$, a secondary index ($\alpha$) is associated
with the multiindex of the geometry tensor $G_K$, and thus the secondary indices indicate
along which dimensions to compute the tensor product $A^0_{i\alpha} G_K^{\alpha}$.
Auxiliary indices ($\beta$) are internal indices within the reference
tensor $A^0$ or the geometry tensor $G_K$ and must be repeated
exactly twice; summation is performed over each auxiliary index
$\beta$ before the tensor product is computed by summation over
secondary indices $\alpha$. Finally, a fixed index is a given constant index
that cannot be evaluated. Fixed indices are used to
represent for example a derivative in a fixed coordinate direction.

Implicit in our algebra is a \emph{grammar} for multilinear forms.  We
could explicitly write an EBNF grammar and use tools such as \texttt{lex}
and \texttt{yacc} to create a compiler for a domain-specific language.
However, by limiting ourselves to overloaded operators, we
successfully embed our language as a high-level library in Python.

To make this concrete, consider test case~2 of section
\ref{sec:examples}, Poisson's equation. The tensor representation
$A^K_i = A^0_{i\alpha} G_K^{\alpha}$ is given by
\begin{equation}
  \begin{split}
    A^0_{i\alpha} &= \int_{K_0}
    \frac{\partial \Phi^1_{i_1}(X)}{\partial X_{\alpha_1}}
    \frac{\partial \Phi^2_{i_2}(X)}{\partial X_{\alpha_2}} \dX, \\
    G_K^{\alpha} &= \det F_K'
    \frac{\partial X_{\alpha_1}}{\partial x_{\beta}}
    \frac{\partial X_{\alpha_2}}{\partial x_{\beta}}.
  \end{split}
\end{equation}
There are here two primary indices ($i_1$ and $i_2$), two secondary
indices ($\alpha_1$ and $\alpha_2$), and one auxiliary index ($\beta$).

\subsection{Evaluation of integrals}

Once the tensor representation (\ref{eq:tensorproduct}) has been
generated, FFC computes all entries of the reference tensor(s) by
quadrature on the reference element. The quadrature rule is
automatically chosen to match the polynomial order of each integrand.
FFC uses FIAT \cite{Kir04,Kir05} as the finite element backend; FIAT
generates the set of basis functions, the quadrature rule, and
evaluates the basis functions and their derivatives at the quadrature
points.

Although FIAT supports many families of finite elements, the current
version of FFC only supports general order continuous/discontinuous
Lagrange finite elements and first order Crouzeix--Raviart finite
elements on triangles and on tetrahedra (or any other finite element
with nodes given by pointwise evaluation). Support for other families
of finite elements will be added in future versions.

Computing integrals is the most expensive step in the compilation of a
form. The typical run-time (of the compiler) ranges between 0.1 and 30 seconds,
depending on the type of form and finite element.

\subsection{Generation of code}

\begin{figure}[htbp]
  \begin{center}
    \includegraphics[width=11cm]{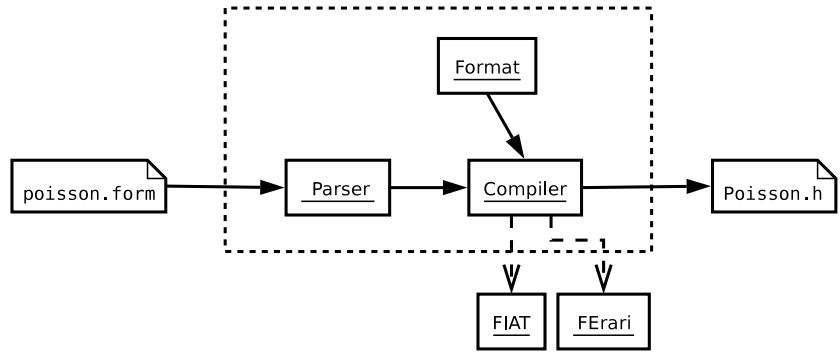}
    \caption{Diagram of the components of the form compiler FFC.}
    \label{fig:components}
  \end{center}
\end{figure}

When a form has been parsed, the tensor representation has been
generated, and all integrals computed, code is generated for
evaluation of geometry tensors and tensor products.  The form
compiler FFC has been designed to allow for generation of code in
multiple different languages. Code is generated according to a
specific \emph{format} (which is essentially a Python dictionary) that
controls the output code being generated, see Figure
\ref{fig:components}. The current version of FFC supports four output
formats: C++ code for DOLFIN~\cite{logg:www:01,logg:preprint:06},
\LaTeX{} code (for verification and presentation purposes), a raw
format that just lists the values of the reference tensors, and the
recently added ASE~format~\cite{www:ASE} for compilation of forms for the
next generation of PETSc. FFC can be easily extended with new output
formats, including for example Python, C, or Fortran.

\subsection{Input/output}

FFC can be used either as a Python package or from the command-line.
We here give a brief description of how FFC can be called from the
command-line to generate C++ code for DOLFIN.
To use FFC from the command-line, one specifies the form in a text
file in a special language for variational forms, which is simply
Python equipped with the hierarchy of classes and operators
discussed above in Section \ref{sec:parsing}. In Table
\ref{tab:poisson,ffc} we give the complete code for the specification
of test case~2, Poisson's equation. Note that FFC uses
tensor-notation, and thus the summation over the index \texttt{i} is
implicit. Also note that the integral over an element $K$ is denoted by
\texttt{*dx}.

\begin{table}[htbp]
  \begin{center}
    \small
    \begin{verbatim}
    element = FiniteElement("Lagrange", "tetrahedron", 3)

    v = BasisFunction(element)
    u = BasisFunction(element)
    f = Function(element)
    i = Index()

    a = v.dx(i)*u.dx(i)*dx
    L = v*f*dx
    \end{verbatim}
    \normalsize
    \caption{The complete code for specification of test case~2,
      Poisson's equation, with piecewise cubics on tetrahedra
      in the language of FFC. Alternatively, the form can be specified
      in terms of standard operators: \texttt{a = dot(grad(v), grad(u))*dx}.}
    \label{tab:poisson,ffc}
  \end{center}
\end{table}

Assuming that the form has been specified in the file
\texttt{Poisson.form}, the form can be compiled using the command
\texttt{ffc Poisson.form}.  This generates the C++ file
\texttt{Poisson.h} to be included in a DOLFIN program.  In Table
\ref{tab:poisson,dolfin}, we include part of the output generated by
FFC with input given by the code from Table \ref{tab:poisson,ffc}. In
addition to this code, FFC generates code for the mapping
$\iota(\cdot, \cdot)$ from local to global degrees of
freedom for each finite element space associated with the form.
Note that the values of the $10 \times 10$ element tensor
$A^K$ (in the case of cubics on triangles) are stored as one contiguous
array (\texttt{block}), since this is what the linear algebra backend
of DOLFIN (PETSc) requires for assembly.

\begin{table}[htbp]
  \begin{center}
    \small
    \begin{verbatim}
void eval(double block[], const AffineMap& map) const
{
  // Compute geometry tensors
  double G0_0_0 = map.det*(map.g00*map.g00 + map.g01*map.g01);
  double G0_0_1 = map.det*(map.g00*map.g10 + map.g01*map.g11);
  double G0_1_0 = map.det*(map.g10*map.g00 + map.g11*map.g01);
  double G0_1_1 = map.det*(map.g10*map.g10 + map.g11*map.g11);

  // Compute element tensor
  block[0]  =  4.249999999999996e-01*G0_0_0 + 4.249999999999995e-01*G0_0_1 +
               4.249999999999995e-01*G0_1_0 + 4.249999999999995e-01*G0_1_1;
  block[1]  = -8.749999999999993e-02*G0_0_0 - 8.749999999999995e-02*G0_0_1;
  block[2]  = -8.750000000000005e-02*G0_1_0 - 8.750000000000013e-02*G0_1_1;
  ...
  block[99] =  4.049999999999997e+00*G0_0_0 + 2.024999999999998e+00*G0_0_1 +
               2.024999999999998e+00*G0_1_0 + 4.049999999999995e+00*G0_1_1;
}
    \end{verbatim}
    \normalsize
    \caption{Part of the code generated by FFC for the input code from
      Table \ref{tab:poisson,ffc}.}
    \label{tab:poisson,dolfin}
  \end{center}
\end{table}

\subsection{Completing the toolchain}
\label{sec:toolchain}

With the FEniCS project \cite{logg:www:03}, we have the beginnings of
a working system realizing (in part) the Automation of Computational
Mathematical Modeling, and the form compiler FFC is just one of several
components needed to complete the toolchain. FIAT automates the
generation of finite element spaces and FFC automates the evaluation
of variational forms. Furthermore, PETSc
\cite{www:PETSc,BalBus04,BalEij97}, automating the solution of
discrete systems, is used as the solver backend of FEniCS. A common
C++ interface to the different FEniCS components is provided by
DOLFIN.

A complete automation of CMM, as outlined in \cite{logg:thesis:03}, is a major
task and we hope that by a modular approach we can contribute to this
automation.

\section{Benchmark results}
\label{sec:benchmarks}

As noted above, the speedup for the code generated by the form
compiler FFC can in many cases be significant. Below, we present a
comparison with the standard quadrature-based approach for the test
cases discussed above in Section \ref{sec:examples}.

The forms were compiled for a range of polynomial degrees using FFC
version~0.1.6. This version of FFC does not take into account any of
the optimizations discussed in Section \ref{sec:optimization}, other
than not generating code for multiplication with zero entries of the
reference tensor.

For the quadrature-based code, all basis functions and their
derivatives were pretabulated at the quadrature points using FIAT.
Loops for all scalar products were completely unrolled.

In all cases, we have used the "collapsed-coordinate" Gauss-Jacobi
rules described by Karniadakis and Sherwin~\cite{KarShe99}.  These
take tensor-product Gaussian integration rules over the square and
cube and map them to the reference simplex.  These rules are not the
best known (see for example~\cite{Dun85}), but they are fairly
simple to generate for arbitrary degree.  Eventually, these rules
will be integrated with FIAT, but even if we reduce the number of
quadrature points by a factor of five, FFC still outperforms
quadrature.

The codes were compiled with gcc (g++) version 3.3.6 and the benchmark
results presented below were obtained on an Intel Pentium 4 (CPU
3.0~GHz, 2GB RAM) running Debian GNU/Linux. The times reported are for
the computation of each entry of the element tensor on one million
elements (scaled). The total time can be obtained by multiplying with
$n^2$, the number of entries of the element tensor.  The complete
source-code for the benchmarks can be obtained from the FEniCS web
site \cite{logg:www:03}.

\subsection{Summary of results}

In Table \ref{tab:speedup}, we summarize the results for test cases
1--4. In all cases, the speedup $T_Q/T_T$ is significant, ranging
between a factor 10--1500.

From Section \ref{sec:complexity}, we know that the speedup for the
mass matrix should grow as $q^d$, but from Table \ref{tab:speedup} it
is clear that the speedup is not quadratic for $d = 2$ and for $d = 3$,
an optimum seems to be reached around $q = 8$.

\begin{figure}[htbp]
  \begin{center}
    \psfrag{q}{$q$}
    \includegraphics[width=11cm]{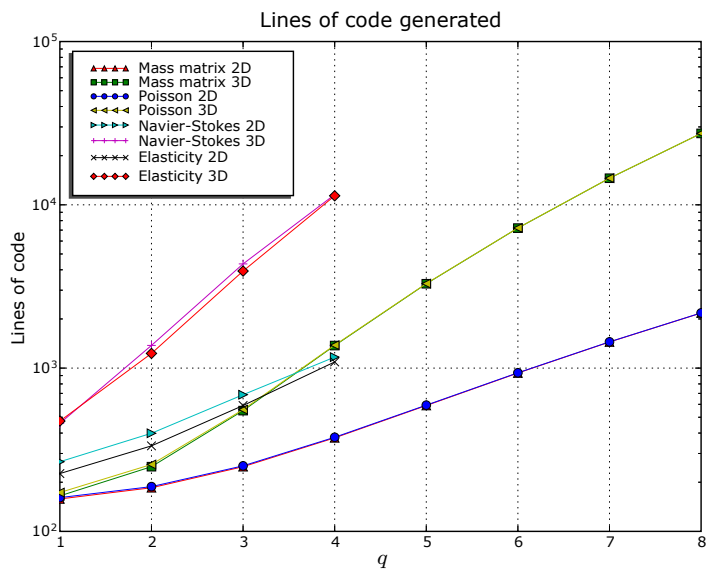}
    \caption{Lines of code generated by the form compiler FFC as
      function of the polynomial degree $q$.}
    \label{fig:klocs}
  \end{center}
\end{figure}

The reason that the predicted speedups are not obtained in practice is
that the complexity estimates presented in Section
\ref{sec:complexity} only account for the number of floating-point
operations. When the polynomial degree $q$ grows, the number of lines
of code generated by the form compiler grows. FFC unrolls
all loops and generates one line of code for each entry of
the element tensor to be computed. For a bilinear form, the number of
entries is $n^2 \sim q^{2d}$. With $q = 8$, the number of lines of
code generated is about $25,000$ for the mass matrix and Poisson in
3D, see Figure \ref{fig:klocs}. As the number of lines of code grows,
memory access becomes more important and dominates the run-time.
Using BLAS to compute tensor products as discussed above might
lead to more efficient memory traffic.

Note however that although the optimal speedup is not obtained, the
speedup is in all cases significant, even at $q = 1$.

\begin{table}[htbp]
  \begin{center}
    \small
    \begin{tabular}{|l|r|r|r|r|r|r|r|r|}
      \hline
      Form & $q = 1$ & $q = 2$ & $q = 3$ & $q = 4$ & $q = 5$ & $q = 6$ & $q = 7$ &$q = 8$ \\
      \hline
      \hline
      Mass 2D           & 12 & 31  & 50  & 78  & 108 & 147 & 183  & 232  \\
      Mass 3D           & 21 & 81  & 189 & 355 & 616 & 881 & 1442 & 1475 \\
      Poisson 2D        & 8  & 29  & 56  & 86  & 129 & 144 & 189  & 236  \\
      Poisson 3D        & 9  & 56  & 143 & 259 & 427 & 341 & 285  & 356  \\
      Navier--Stokes 2D & 32 & 33  & 53  & 37  & --- & --- & --- & ---   \\
      Navier--Stokes 3D & 77 & 100 & 61  & 42  & --- & --- & --- & ---   \\
      Elasticity 2D     & 10 & 43  & 67  & 97  & --- & --- & --- & ---   \\
      Elasticity 3D     & 14 & 87  & 103 & 134 & --- & --- & --- & ---   \\
      \hline
    \end{tabular}
    \normalsize
    \caption{Speedups $T_Q/T_T$ for test cases 1--4 in 2D and 3D.}
    \label{tab:speedup}
  \end{center}
\end{table}

\subsection{Results for test cases}

In Figures \ref{fig:result,1}--\ref{fig:result,4}, we present the
results for test cases~1--4 discussed above in Section
\ref{sec:examples}. In connection to each of the results, we include
the specification of the form in the language used by the form
compiler FFC. Because of limitations in the current implementation of
FFC, the comparison is made for polynomial order $q \leq 8$ in test
cases~1--2 and $q \leq 4$ in test cases~3--4. Higher polynomial order
is possible but is very costly to compile (compare Figure~\ref{fig:klocs}).

\begin{figure}[htbp]
  \begin{center}
    \psfrag{q}{$q$}
    \includegraphics[width=11cm]{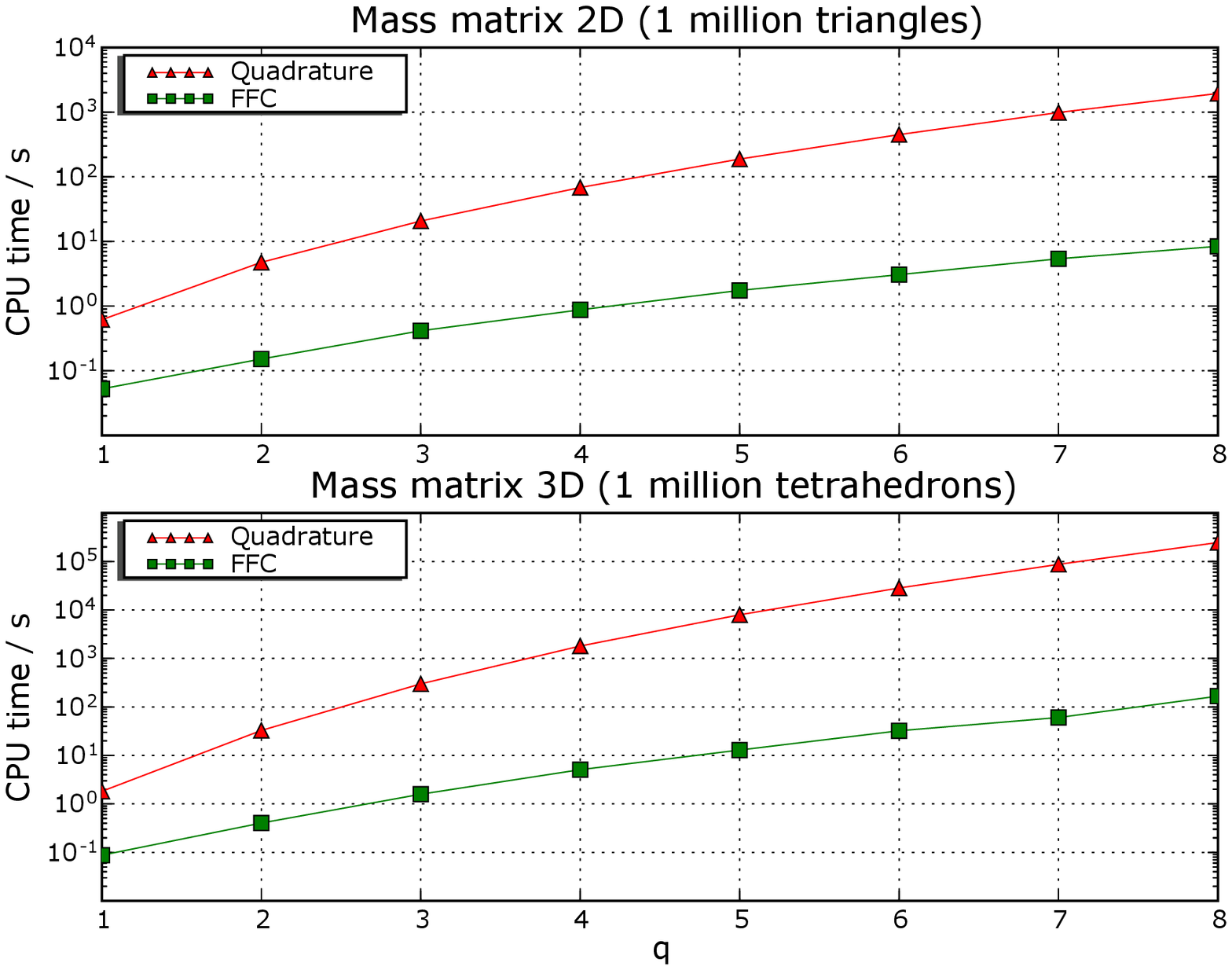}
    \caption{Benchmark results for test case~1, the mass matrix,
      specified in FFC by \texttt{a = v*u*dx}.}
    \label{fig:result,1}
  \end{center}
\end{figure}

\begin{figure}[htbp]
  \begin{center}
    \psfrag{q}{$q$}
    \includegraphics[width=11cm]{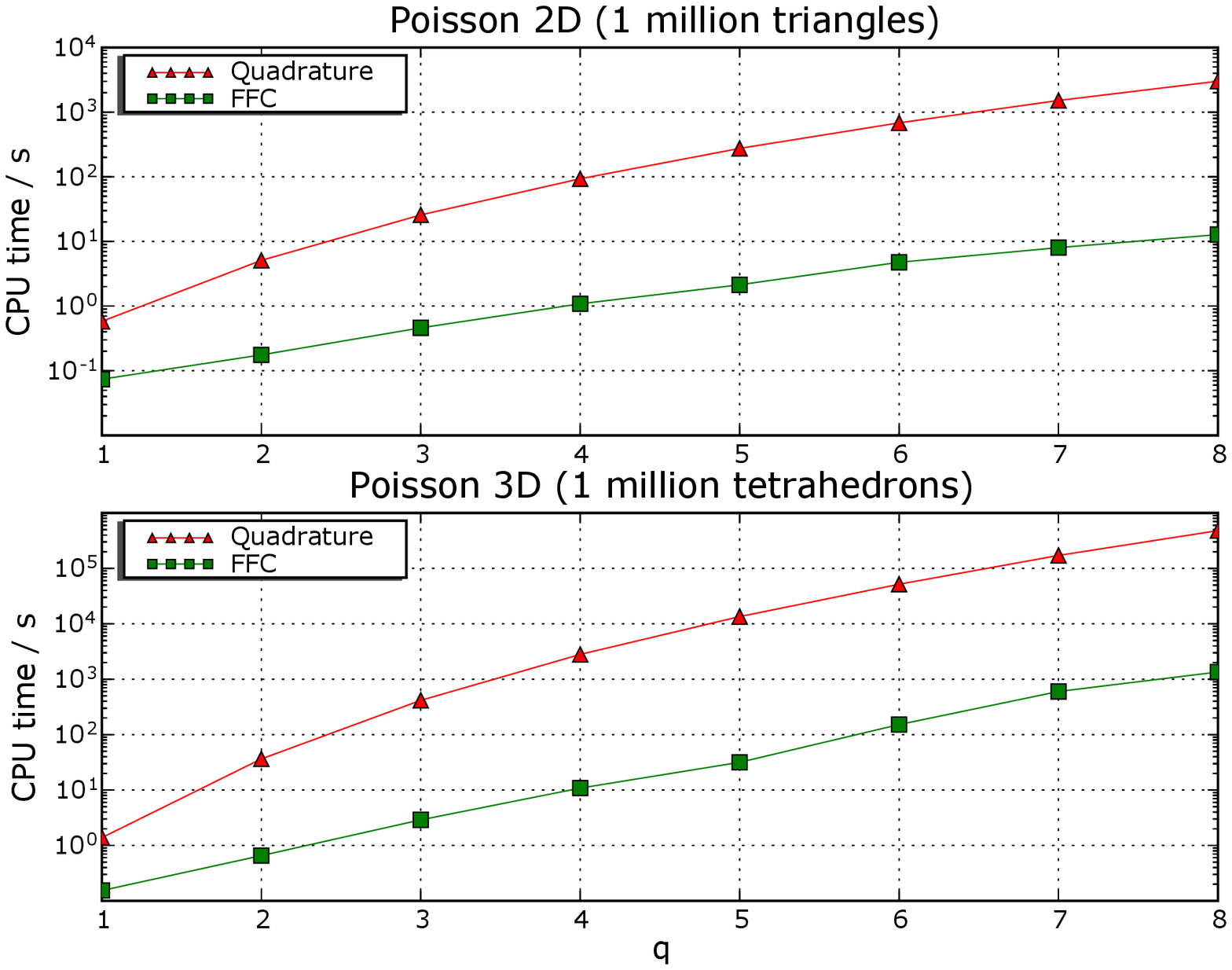}
    \caption{Benchmark results for test case~2, Poisson's equation,
      specified in FFC by \texttt{a = v.dx(i)*u.dx(i)*dx}.}
    \label{fig:result,2}
  \end{center}
\end{figure}

\begin{figure}[htbp]
  \begin{center}
    \psfrag{q}{$q$}
    \includegraphics[width=11cm]{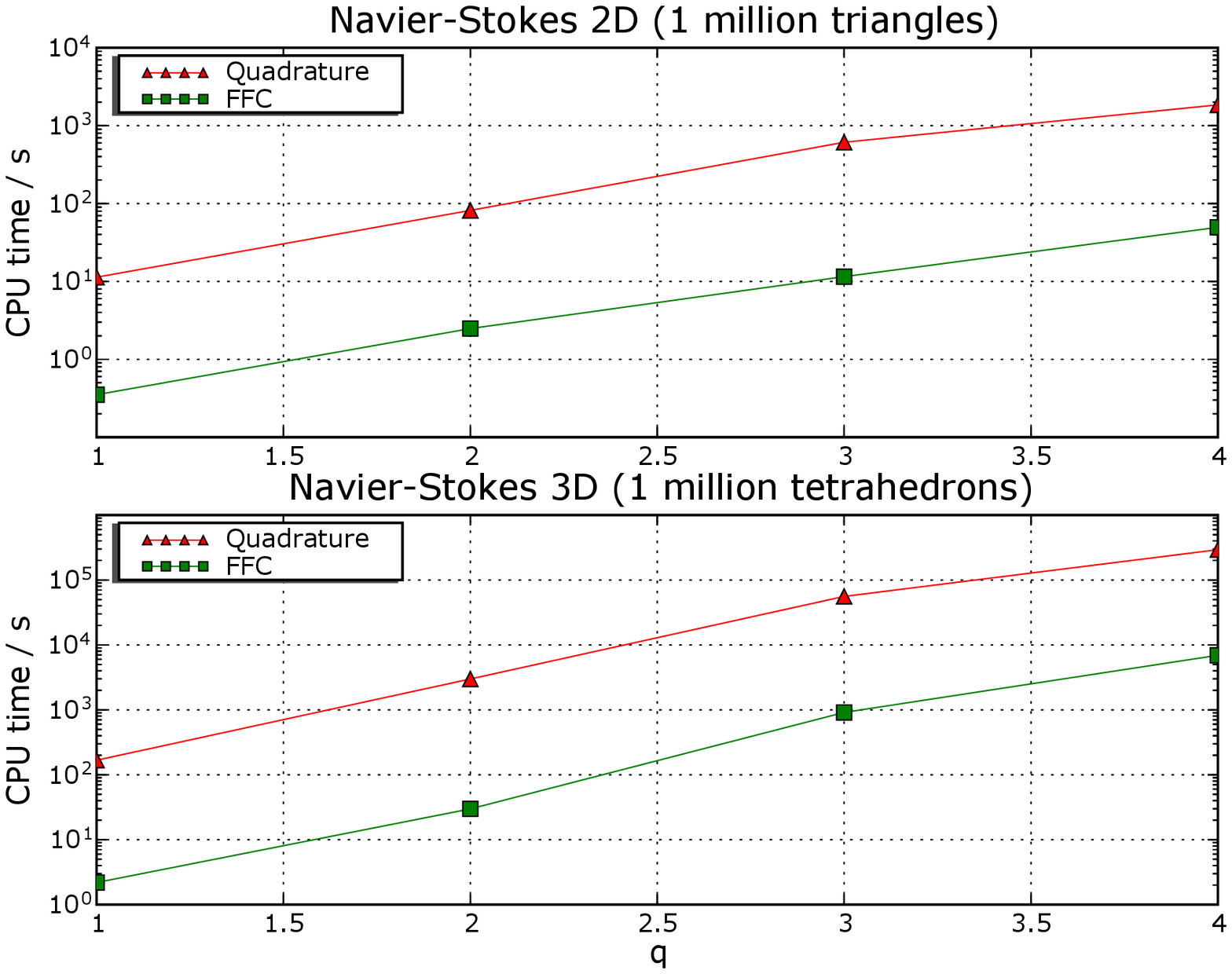}
    \caption{Benchmark results for test case~3, the nonlinear term of
      the incompressible Navier--Stokes equations, specified in FFC
      by \texttt{a = v[i]*w[j]*u[i].dx(j)*dx}.}
    \label{fig:result,3}
  \end{center}
\end{figure}

\begin{figure}[htbp]
  \begin{center}
    \psfrag{q}{$q$}
    \includegraphics[width=11cm]{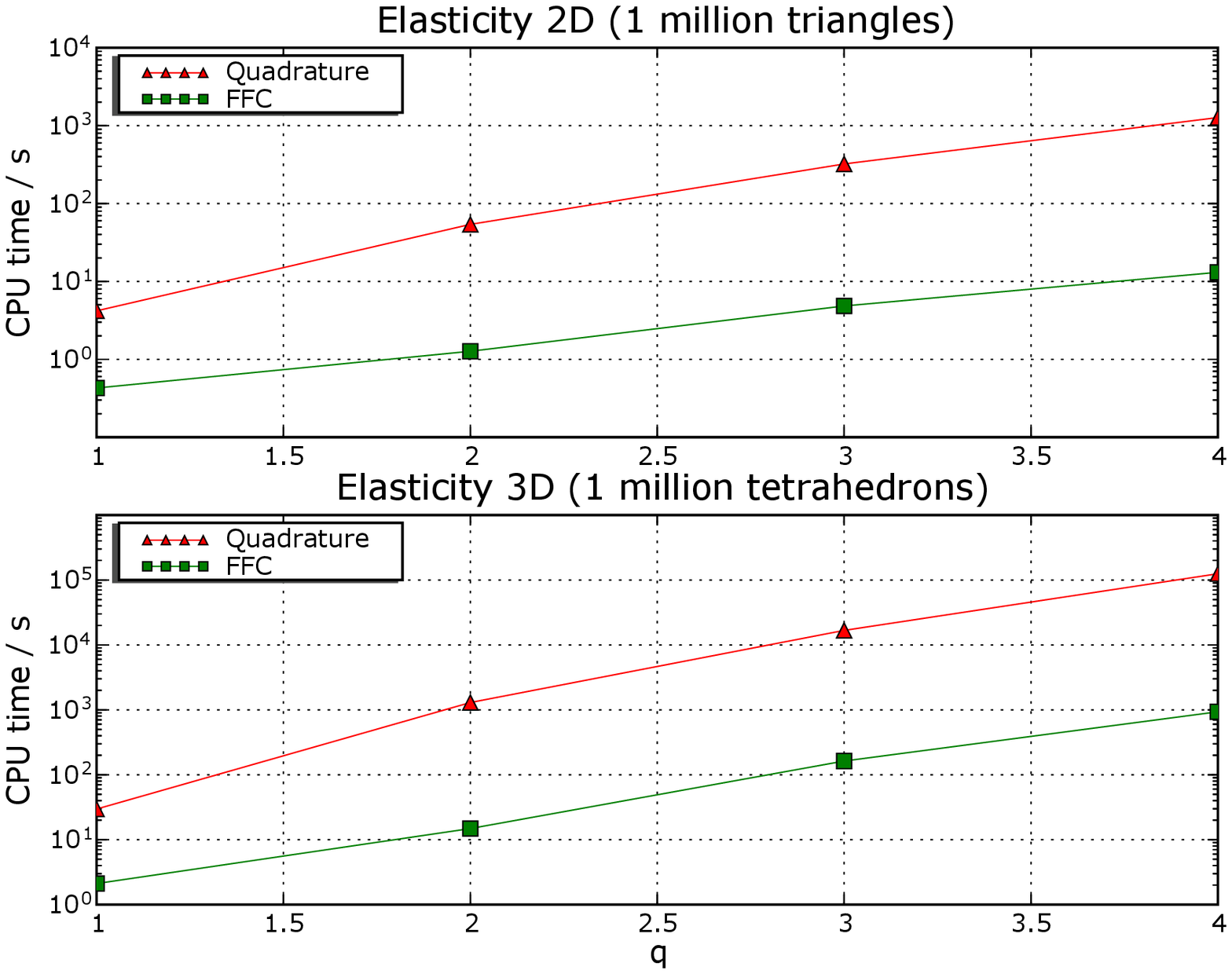}
    \caption{Benchmark results for test case~4, the strain-strain term
      of linear elasticity, specified in FFC
      by \texttt{a = 0.25*(v[i].dx(j) + v[j].dx(i)) * (u[i].dx(j) + u[j].dx(i)) * dx}.}
    \label{fig:result,4}
  \end{center}
\end{figure}

\section{Concluding remarks and future directions}
\label{sec:conclusion}

We have demonstrated a proof-of-concept form compiler that for a wide
range of variational forms can generate code that gives significant
speedups compared to the standard quadrature-based approach.

The form compiler FFC is still in its early stages of development but
is already in production use. A number of basic modules based on FFC
have been implemented in DOLFIN and others are currently being
developed (Navier--Stokes and updated elasticity). This will serve as
a test bed for future development of FFC.

Future plans for FFC include adding support for integrals over the
boundary (adding the operator \texttt{*ds} to the language), support
for automatic differentiation of nonlinear forms and automatic
generation of dual problems and a posteriori error estimators
\cite{EriEst95,BecRan01}, optimization through FErari
\cite[SISC]{logg:article:07}, \cite[BIT]{KirKne05}, adding support for
new families of finite elements, including elements that require
non-affine mappings from the reference element. In addition to general
order continuous/discontinuous Lagrange finite elements and
Crouzeix--Raviart \cite{CroRav73} finite elements, the plan is to add
support for Raviart--Thomas \cite{RavTho77a}, Nedelec \cite{Ned80},
Brezzi--Douglas--Marini \cite{BreDou85},
Brezzi--Douglas--Fortin--Marini \cite{BreFor91},
Taylor--Hood \cite{Bof97,BreSco94}, and
Arnold--Winther \cite{ArnWin02} elements.

Furthermore, the fact that Python is an interpreted language does
impose a penalty on the performance of the compiler (but not on the
generated code). However, this can be overcome by more extensive use
of BLAS in numerically intensive parts of the compiler, such as the
precomputation of integrals on the reference element.
We also plan to investigate the use of BLAS for evaluation of tensor
products as an alternative to generating explicit unrolled code.
Other topics of interest include automatic verification of the
correctness of the code generated by the form compiler \cite[Verification]{KirStr05}.

\begin{ack}

We wish to thank the FEniCS team, in particular Johan Hoffman, Johan
Jansson, Claes Johnson, Matthew Knepley, and Ridgway Scott, for
substantial suggestions and comments regarding this paper.  We also
want to thank one of the referees for pointing out how the tensor
representation can be generalized to a quadrature-based evaluation
scheme for non-affine mappings.

\end{ack}
\bibliographystyle{acmtrans}
\bibliography{bibliography}

\begin{thebibliography}{}

\bibitem[\protect\citeauthoryear{Arnold and Winther}{Arnold and
  Winther}{2002}]{ArnWin02}
{\sc Arnold, D.~N.} {\sc and} {\sc Winther, R.} 2002.
\newblock Mixed finite elements for elasticity.
\newblock {\em Numer. Math.\/}~{\em 92,\/}~3, 401--419.

\bibitem[\protect\citeauthoryear{Bagheri and Scott}{Bagheri and
  Scott}{2003}]{www:Analysa}
{\sc Bagheri, B.} {\sc and} {\sc Scott, R.} 2003.
\newblock Analysa.
\newblock \url{http://people.cs.uchicago.edu/~ridg/al/aa.html}.

\bibitem[\protect\citeauthoryear{Balay, Buschelman, Eijkhout, Gropp, Kaushik,
  Knepley, McInnes, Smith, and Zhang}{Balay et~al\mbox{.}}{2004}]{BalBus04}
{\sc Balay, S.}, {\sc Buschelman, K.}, {\sc Eijkhout, V.}, {\sc Gropp, W.~D.},
  {\sc Kaushik, D.}, {\sc Knepley, M.~G.}, {\sc McInnes, L.~C.}, {\sc Smith,
  B.~F.}, {\sc and} {\sc Zhang, H.} 2004.
\newblock {PETS}c users manual.
\newblock Tech. Rep. ANL-95/11 - Revision 2.1.5, Argonne National Laboratory.

\bibitem[\protect\citeauthoryear{Balay, Buschelman, Gropp, Kaushik, Knepley,
  McInnes, Smith, and Zhang}{Balay et~al\mbox{.}}{2005a}]{www:ASE}
{\sc Balay, S.}, {\sc Buschelman, K.}, {\sc Gropp, W.~D.}, {\sc Kaushik, D.},
  {\sc Knepley, M.~G.}, {\sc McInnes, L.~C.}, {\sc Smith, B.~F.}, {\sc and}
  {\sc Zhang, H.} 2005a.
\newblock Anl anl sidl environmentsidl environment.
\newblock \url{http://www-unix.mcs.anl.gov/ase/}.

\bibitem[\protect\citeauthoryear{Balay, Buschelman, Gropp, Kaushik, Knepley,
  McInnes, Smith, and Zhang}{Balay et~al\mbox{.}}{2005b}]{www:PETSc}
{\sc Balay, S.}, {\sc Buschelman, K.}, {\sc Gropp, W.~D.}, {\sc Kaushik, D.},
  {\sc Knepley, M.~G.}, {\sc McInnes, L.~C.}, {\sc Smith, B.~F.}, {\sc and}
  {\sc Zhang, H.} 2005b.
\newblock {PETS}c.
\newblock \url{http://www.mcs.anl.gov/petsc/}.

\bibitem[\protect\citeauthoryear{Balay, Eijkhout, Gropp, McInnes, and
  Smith}{Balay et~al\mbox{.}}{1997}]{BalEij97}
{\sc Balay, S.}, {\sc Eijkhout, V.}, {\sc Gropp, W.~D.}, {\sc McInnes, L.~C.},
  {\sc and} {\sc Smith, B.~F.} 1997.
\newblock Efficient management of parallelism in object oriented numerical
  software libraries.
\newblock In {\em Modern Software Tools in Scientific Computing}, {E.~Arge},
  {A.~M. Bruaset}, {and} {H.~P. Langtangen}, Eds. Birkh{\"{a}}user Press,
  163--202.

\bibitem[\protect\citeauthoryear{Bangerth, Hartmann, and Kanschat}{Bangerth
  et~al\mbox{.}}{2005}]{www:deal}
{\sc Bangerth, W.}, {\sc Hartmann, R.}, {\sc and} {\sc Kanschat, G.} 2005.
\newblock {\tt deal.{I}{I}} {D}ifferential {E}quations {A}nalysis {L}ibrary.
\newblock \url{http://www.dealii.org}.

\bibitem[\protect\citeauthoryear{Becker and Rannacher}{Becker and
  Rannacher}{2001}]{BecRan01}
{\sc Becker, R.} {\sc and} {\sc Rannacher, R.} 2001.
\newblock An optimal control approach to a posteriori error estimation in
  finite element methods.
\newblock {\em Acta Numerica\/}~{\em 10}, 1--102.

\bibitem[\protect\citeauthoryear{Boffi}{Boffi}{1997}]{Bof97}
{\sc Boffi, D.} 1997.
\newblock Three-dimensional finite element methods for the {S}tokes problem.
\newblock {\em SIAM J. Numer. Anal.\/}~{\em 34,\/}~2, 664--670.

\bibitem[\protect\citeauthoryear{Brenner and Scott}{Brenner and
  Scott}{1994}]{BreSco94}
{\sc Brenner, S.~C.} {\sc and} {\sc Scott, L.~R.} 1994.
\newblock {\em The Mathematical Theory of Finite Element Methods}.
\newblock Springer-Verlag.

\bibitem[\protect\citeauthoryear{Brezzi, Douglas, and Marini}{Brezzi
  et~al\mbox{.}}{1985}]{BreDou85}
{\sc Brezzi, F.}, {\sc Douglas, Jr., J.}, {\sc and} {\sc Marini, L.~D.} 1985.
\newblock Two families of mixed finite elements for second order elliptic
  problems.
\newblock {\em Numer. Math.\/}~{\em 47,\/}~2, 217--235.

\bibitem[\protect\citeauthoryear{Brezzi and Fortin}{Brezzi and
  Fortin}{1991}]{BreFor91}
{\sc Brezzi, F.} {\sc and} {\sc Fortin, M.} 1991.
\newblock {\em Mixed and hybrid finite element methods}. Springer Series in
  Computational Mathematics, vol.~15.
\newblock Springer-Verlag, New York.

\bibitem[\protect\citeauthoryear{Castillo, Koning, Rieben, and White}{Castillo
  et~al\mbox{.}}{2004}]{CasKon04}
{\sc Castillo, P.}, {\sc Koning, J.}, {\sc Rieben, R.}, {\sc and} {\sc White,
  D.} 2004.
\newblock A discrete differential forms framework for computational
  electromagnetics.
\newblock {\em Computer Modeling in Engineering and Sciences\/}~{\em 5,\/}~4,
  331 -- 346.

\bibitem[\protect\citeauthoryear{Castillo, Rieben, and White}{Castillo
  et~al\mbox{.}}{2005}]{CasRie05}
{\sc Castillo, P.}, {\sc Rieben, R.}, {\sc and} {\sc White, D.} 2005.
\newblock Femster: an object-oriented class library of discrete differential
  forms.
\newblock {\em To appear in ACM Trans. Math. Software\/}.

\bibitem[\protect\citeauthoryear{Ciarlet}{Ciarlet}{1976}]{Cia76}
{\sc Ciarlet, P.~G.} 1976.
\newblock {\em Numerical Analysis of the Finite Element Method}.
\newblock Les Presses de l'Universite de Montreal.

\bibitem[\protect\citeauthoryear{Crouzeix and Raviart}{Crouzeix and
  Raviart}{1973}]{CroRav73}
{\sc Crouzeix, M.} {\sc and} {\sc Raviart, P.~A.} 1973.
\newblock Conforming and nonconforming finite element methods for solving the
  stationary stokes equations.
\newblock {\em RAIRO Anal. Numér.\/}~{\em 7}, 33--76.

\bibitem[\protect\citeauthoryear{Dular and Geuzaine}{Dular and
  Geuzaine}{2005}]{www:GetDP}
{\sc Dular, P.} {\sc and} {\sc Geuzaine, C.} 2005.
\newblock Get{DP}: a {G}eneral environment for the treatment of {D}iscrete
  {P}roblems.
\newblock \url{http://www.geuz.org/getdp/}.

\bibitem[\protect\citeauthoryear{Dunavant}{Dunavant}{1985}]{Dun85}
{\sc Dunavant, D.~A.} 1985.
\newblock High degree efficient symmetrical {G}aussian quadrature rules for the
  triangle.
\newblock {\em Internat. J. Numer. Methods Engrg.\/}~{\em 21,\/}~6, 1129--1148.

\bibitem[\protect\citeauthoryear{Eriksson, Estep, Hansbo, and Johnson}{Eriksson
  et~al\mbox{.}}{1995}]{EriEst95}
{\sc Eriksson, K.}, {\sc Estep, D.}, {\sc Hansbo, P.}, {\sc and} {\sc Johnson,
  C.} 1995.
\newblock Introduction to adaptive methods for differential equations.
\newblock {\em Acta Numerica\/}~{\em 4}, 105--158.

\bibitem[\protect\citeauthoryear{Eriksson, Estep, Hansbo, and Johnson}{Eriksson
  et~al\mbox{.}}{1996}]{EriEst96}
{\sc Eriksson, K.}, {\sc Estep, D.}, {\sc Hansbo, P.}, {\sc and} {\sc Johnson,
  C.} 1996.
\newblock {\em Computational Differential Equations}.
\newblock Cambridge University Press.

\bibitem[\protect\citeauthoryear{Eriksson, Estep, and Johnson}{Eriksson
  et~al\mbox{.}}{2003}]{EriEst03c}
{\sc Eriksson, K.}, {\sc Estep, D.}, {\sc and} {\sc Johnson, C.} 2003.
\newblock {\em Applied Mathematics: Body and Soul}. Vol. III.
\newblock Springer-Verlag.

\bibitem[\protect\citeauthoryear{{F}ree~{S}oftware
  {F}oundation}{{F}ree~{S}oftware {F}oundation}{1991}]{www:GPL}
{\sc {F}ree~{S}oftware {F}oundation}. 1991.
\newblock {GNU} {GPL}.
\newblock \url{http://www.gnu.org/copyleft/gpl.html}.

\bibitem[\protect\citeauthoryear{Hoffman, Jansson, Johnson, Knepley, Kirby,
  Logg, and Scott}{Hoffman et~al\mbox{.}}{2005}]{logg:www:03}
{\sc Hoffman, J.}, {\sc Jansson, J.}, {\sc Johnson, C.}, {\sc Knepley, M.},
  {\sc Kirby, R.~C.}, {\sc Logg, A.}, {\sc and} {\sc Scott, L.~R.} 2005.
\newblock {\em {FE}ni{CS}}.
\newblock \url{http://www.fenics.org/}.

\bibitem[\protect\citeauthoryear{Hoffman, Jansson, and Logg}{Hoffman
  et~al\mbox{.}}{2005}]{logg:www:01}
{\sc Hoffman, J.}, {\sc Jansson, J.}, {\sc and} {\sc Logg, A.} 2005.
\newblock {\em {DOLFIN}}.
\newblock \url{http://www.fenics.org/dolfin/}.

\bibitem[\protect\citeauthoryear{Hoffman and Logg}{Hoffman and
  Logg}{2002}]{logg:preprint:06}
{\sc Hoffman, J.} {\sc and} {\sc Logg, A.} 2002.
\newblock {DOLFIN}: {D}ynamic {O}bject oriented {L}ibrary for {FIN}ite element
  computation.
\newblock Tech. Rep. 2002--06, Chalmers Finite Element Center Preprint Series.

\bibitem[\protect\citeauthoryear{Hughes}{Hughes}{1987}]{Hug87}
{\sc Hughes, T. J.~R.} 1987.
\newblock {\em The Finite Element Method: Linear Static and Dynamic Finite
  Element Analysis}.
\newblock Prentice-Hall.

\bibitem[\protect\citeauthoryear{Karniadakis and Sherwin}{Karniadakis and
  Sherwin}{1999}]{KarShe99}
{\sc Karniadakis, G.~E.} {\sc and} {\sc Sherwin, S.~J.} 1999.
\newblock {\em Spectral/{$hp$} element methods for {CFD}}.
\newblock Numerical Mathematics and Scientific Computation. Oxford University
  Press, New York.

\bibitem[\protect\citeauthoryear{Kirby}{Kirby}{2004}]{Kir04}
{\sc Kirby, R.~C.} 2004.
\newblock {FIAT}: A new paradigm for computing finite element basis functions.
\newblock {\em ACM Trans. Math. Software\/}~{\em 30}, 502--516.

\bibitem[\protect\citeauthoryear{Kirby}{Kirby}{2005}]{Kir05}
{\sc Kirby, R.~C.} 2005.
\newblock Optimizing {FIAT} with the level 3 {BLAS}.
\newblock {\em submitted to ACM Trans. Math. Software\/}.

\bibitem[\protect\citeauthoryear{Kirby, Knepley, Logg, and Scott}{Kirby
  et~al\mbox{.}}{2005}]{logg:article:07}
{\sc Kirby, R.~C.}, {\sc Knepley, M.}, {\sc Logg, A.}, {\sc and} {\sc Scott,
  L.~R.} 2005.
\newblock Optimizing the evaluation of finite element matrices.
\newblock {\em To appear in SIAM J. Sci. Comput.\/}.

\bibitem[\protect\citeauthoryear{Kirby, Knepley, and Scott}{Kirby
  et~al\mbox{.}}{2005}]{KirKne05}
{\sc Kirby, R.~C.}, {\sc Knepley, M.}, {\sc and} {\sc Scott, L.~R.} 2005.
\newblock Evaluation of the action of finite element operators.
\newblock {\em submitted to BIT\/}.

\bibitem[\protect\citeauthoryear{Kirby, Strout, Hovland, and Scott}{Kirby
  et~al\mbox{.}}{2005}]{KirStr05}
{\sc Kirby, R.~C.}, {\sc Strout, M.~M.}, {\sc Hovland, P.}, {\sc and} {\sc
  Scott, L.~R.} 2005.
\newblock Verification of scientific code using rationality analysis.
\newblock {\em in preparation\/}.

\bibitem[\protect\citeauthoryear{Langtangen}{Langtangen}{1999}]{Lan99}
{\sc Langtangen, H.~P.} 1999.
\newblock {\em Computational Partial Differential Equations -- Numerical
  Methods and Diffpack Programming}.
\newblock Lecture Notes in Computational Science and Engineering. Springer.

\bibitem[\protect\citeauthoryear{Logg}{Logg}{2004}]{logg:thesis:03}
{\sc Logg, A.} 2004.
\newblock Automation of computational mathematical modeling.
\newblock Ph.D. thesis, Chalmers University of Technology, Sweden.

\bibitem[\protect\citeauthoryear{Long}{Long}{2003}]{Lon03}
{\sc Long, K.} 2003.
\newblock Sundance, a rapid prototyping tool for parallel {PDE}-constrained
  optimization.
\newblock In {\em Large-Scale PDE-Constrained Optimization}. Lecture notes in
  computational science and engineering. Springer-Verlag.

\bibitem[\protect\citeauthoryear{Mackie}{Mackie}{1992}]{Mac92}
{\sc Mackie, R.~I.} 1992.
\newblock Object oriented programming of the finite element method.
\newblock {\em Int. J. Num. Meth. Eng.\/}~{\em 35,\/}~425--436.

\bibitem[\protect\citeauthoryear{Masters, Usmani, Cross, and Lewis}{Masters
  et~al\mbox{.}}{1997}]{MasUsm97}
{\sc Masters, I.}, {\sc Usmani, A.~S.}, {\sc Cross, J.~T.}, {\sc and} {\sc
  Lewis, R.~W.} 1997.
\newblock Finite element analysis of solidification using obect-oriented and
  parallel techniques.
\newblock {\em Int. J. Numer. Meth. Eng.\/}~{\em 40}, 2891--2909.

\bibitem[\protect\citeauthoryear{N{\'e}d{\'e}lec}{N{\'e}d{\'e}lec}{1980}]{Ned8%
0}
{\sc N{\'e}d{\'e}lec, J.-C.} 1980.
\newblock Mixed finite elements in {${\bf R}\sp{3}$}.
\newblock {\em Numer. Math.\/}~{\em 35,\/}~3, 315--341.

\bibitem[\protect\citeauthoryear{Pironneau, Hecht, and Hyaric}{Pironneau
  et~al\mbox{.}}{2005}]{www:FreeFEM}
{\sc Pironneau, O.}, {\sc Hecht, F.}, {\sc and} {\sc Hyaric, A.~L.} 2005.
\newblock Free{FEM}.
\newblock \url{http://www.freefem.org/}.

\bibitem[\protect\citeauthoryear{Raviart and Thomas}{Raviart and
  Thomas}{1977}]{RavTho77a}
{\sc Raviart, P.-A.} {\sc and} {\sc Thomas, J.~M.} 1977.
\newblock A mixed finite element method for 2nd order elliptic problems.
\newblock In {\em Mathematical aspects of finite element methods (Proc. Conf.,
  Consiglio Naz. delle Ricerche (C.N.R.), Rome, 1975)}. Springer, Berlin,
  292--315. Lecture Notes in Math., Vol. 606.

\end{thebibliography}

\begin{received}
\end{received}

\end{document}